\documentclass[a4paper,11pt]{amsart}

\providecommand{\MR}{\relax\ifhmode\unskip\space\fi MR }

\providecommand{\href}[2]{#2}

\usepackage{enumerate, booktabs, mathrsfs}
\usepackage{amsmath, amsfonts, amssymb, amsthm}
\usepackage[all]{xy}
\usepackage[height=22.5cm, width=15.5cm]{geometry}
\usepackage[active]{srcltx}
\usepackage{color}  
\numberwithin{equation}{section}
\newtheorem{thm}{Theorem}[section]
\newtheorem{lemma}[thm]{Lemma}
\newtheorem{proposition}[thm]{Proposition}
\newtheorem{cor}[thm]{Corollary}
\theoremstyle{definition}
\newtheorem{definition}[thm]{Definition}
\newtheorem{remark}[thm]{Remark}  
  
\numberwithin{equation}{section}

\author{Bruce Gilligan}
\address{Dept. of Mathematics \& Statistics, University of Regina, Regina, Canada S4S 0A2}
\email{Bruce.Gilligan@uregina.ca}

\date{November 18, 2026}   

\begin{document}

\title[The holomorphic reduction of pseudoconvex homogeneous manifolds]
{The holomorphic reduction of pseudoconvex homogeneous manifolds}

\begin{abstract}    
Suppose $G$ is a connected complex Lie group and $H$ is a closed complex subgroup.  
Then there exists a closed complex subgroup $J$ of $G$ containing $H$ such that 
the fibration $\pi:G/H \to G/J$ is the holomorphic reduction of $G/H$, i.e., $G/J$ is holomorphically 
separable and ${\mathcal O}(G/H) \cong \pi^*{\mathcal O}(G/J)$.  
In this paper we prove that if $G/H$ is pseudoconvex, i.e.,    if      
$G/H$ admits a continuous plurisubharmonic exhaustion function,           
then $G/J$ is Stein and $J/H$ has no non-constant holomorphic functions.   
\end{abstract}  

\thanks{This work was partially supported by an NSERC Discovery Grant.
The author used Claude (Anthropic) to assist with reviewing the manuscript, checking and repairing arguments, and editing the exposition. 
All mathematical content was verified by the author, who takes full responsibility for it.    
We thank  A. Huckleberry, C. Miebach, and K. Oeljeklaus for discussions concerning this paper. } 

\maketitle

\section{Introduction}   \label{sect1}  

\medskip 
In this paper we consider complex homogeneous manifolds $X := G/H$, where
$G$ is a connected complex Lie group and $H$ a closed complex subgroup.    
Setting 
\[    
       J \; := \;  \{ \ g \in G \   | \  f(gH) \;  = \;  f(eH) \mbox{  for all  } f\in{\mathcal O}(G/H) \  \}   
\]           
yields a closed complex subgroup of $G$ containing $H$ and consequently one has   
the holomorphic fibration $\pi : G/H \to G/J, \; gH \mapsto gJ$   
which we call the \textbf{holomorphic reduction} of $G/H$   \cite{GH78}.  
By construction $G/J$ is holomorphically separable and ${\mathcal O}(G/H) \cong \pi^* {\mathcal O}(G/J)$.   
     
\medskip 
In general the fiber $J/H$ of the holomorphic reduction may carry non-constant holomorphic functions, and may even be Stein. 
For instance, if $\Gamma\cong\mathbb Z$ is a discrete subgroup of the unipotent subgroup $U$ of upper triangular matrices in $SL(2,\mathbb C)$, then $U/\Gamma \cong \mathbb C^*$ is the fiber of  
the holomorphic reduction of $SL(2,\mathbb C)/\Gamma$ 
which is thus a $\mathbb C^*$-bundle over $\mathbb C^2\setminus\{0\}$ \cite[\S 1.4]{BO73}; 
see \cite[\S 4]{Loe85} for an example involving solvmanifolds. 
Part of the content of Theorem \ref{mnthm} is that this cannot happen when $G/H$ is pseudoconvex.
  
\medskip  
As far as the base of the holomorphic reduction is concerned we remark that 
holomorphic separability implies Steinness for complex Lie groups \cite{MM60}, for complex nilmanifolds \cite{GH78}, 
and for complex solvmanifolds \cite{HO86}.  
But the situation is much different when the group acting is semisimple or reductive.  
For example, $\mathbb C^n \setminus \{ 0 \}$ is not Stein whenever $n>1$.  
What is known is that $G/H$ is holomorphically separable for $G$ reductive implies $H$ is an 
algebraic subgroup of $G$ \cite[\S 2.5]{BO73} and such a $G/H$ is Stein if and only if $H$ is reductive \cite{Mat60} and \cite{On60}.  
A reductive complex Lie group is the complexification of a totally real maximal compact subgroup and the  
compactness of this subgroup plays an essential role.  
An excellent survey of the function theory on $G/H$ for $G$ reductive can be found in \cite[Chapter 5]{Akh95} and we refer   
the interested reader to that book and the references listed there.        

\medskip   
In the present paper we analyze the structure of the holomorphic reductions of \textbf{pseudoconvex} complex homogeneous manifolds,   
where a complex manifold is pseudoconvex if it admits a continuous plurisubharmonic exhaustion function.       
In \cite[Main Theorem]{GMO13} we proved that the base of the holomorphic reduction of any pseudoconvex complex   
homogeneous manifold $G/H$ is Stein and its fiber has no non-constant holomorphic functions if $G$ is solvable  (see also  \cite[Theorem 8.5]{GMO13}) 
or if $G$ is reductive (see also \cite[Theorem 7.4]{GMO13}).         
Here we prove that this result holds more generally for pseudoconvex homogeneous manifolds of mixed groups, where by a mixed group       
we mean that $G$ is a connected, simply connected, complex Lie group that 
has a Levi-Malcev decomposition $G = S \ltimes R$, where $R$ is the radical of $G$,    
$S$ is a maximal semisimple subgroup, and both $\dim R > 0$ and $\dim S > 0$.   
In Section 3 we will see that the general case reduces to the solvable, semisimple and mixed cases. 

\begin{thm} \label{mnthm}  
Let $G$ be a connected complex Lie group and $H$ a closed complex subgroup of $G$ such 
that $G/H$ is pseudoconvex.  
Suppose $G/H \to G/J$ is the holomorphic reduction of $G/H$.  
Then $G/J$ is Stein and ${\mathcal O}(J/H) = \mathbb C$.  
\end{thm}   

Theorem \ref{mnthm} was stated as Theorem 1.1 in our earlier paper \cite{Gil17}, but the proof given there is incomplete.  
Lemma 2.6 of \cite{Gil17} is false as stated, and the proof of Lemma 2.7 there depends on it; 
the argument in \S 3.2 of \cite{Gil17} does not show that the orbits of the complexified radical are closed; 
the case where the identity component of $L$ is semisimple but its complex hull is not was overlooked; 
and the Hirschowitz annihilator has to be defined with a complex parameter, as in \cite{GMO13}.  
The present paper gives a complete and self-contained proof of Theorem \ref{mnthm}, 
which supersedes the argument in \cite{Gil17}. 

\medskip 
In Section 2 we collect the technical results that are needed later. 
Section 3 contains the proof of Theorem \ref{mnthm}, and its strategy is described in \S 3.1.

\section{Technical Preparations}  

\subsection{The Hirschowitz annihilator}   \label{Hirschowitz}  
Every element $\xi\in\mathfrak g$, the Lie algebra of $G$, can be thought of as a right invariant vector field on $G$ 
and, as such, pushes down to a holomorphic vector field $\xi_X$ on any complex homogeneous manifold $X:=G/H$ of the group $G$.     
An inner integral curve in such a homogeneous space $X$ is a non-constant holomorphic map $\mathbb C \to X$ with relatively compact image 
in $X$ that is the integral curve of some vector field $\xi_X$ associated to some $\xi\in\mathfrak g$.  
Hirschowitz \cite{Hir75} considered such concepts in the context of infinitesimally homogeneous manifolds, where a manifold  
is infinitesimally homogeneous if every tangent space is generated by global holomorphic vector fields.  
A complex manifold that is homogeneous under the action of a Lie group of holomorphic transformations is 
infinitesimally homogeneous.     

\medskip   
Hirschowitz showed that if a pseudoconvex, infinitesimally homogeneous $X$ does not contain    
any inner integral curve, then it is Stein \cite[Proposition 3.4]{Hir75}.   
This is the starting point of our investigations of pseudoconvex homogeneous manifolds that are not Stein that was presented in \cite{GMO13}.      
From Liouville's theorem for subharmonic functions bounded above on $\mathbb C$ 
one gets that any plurisubharmonic function on a complex manifold $X$ is constant along every inner integral curve in $X$.   
One has to determine the ``directions of degeneracy'' of plurisubharmonic functions in terms of a certain subset of 
$\mathfrak g$ whose corresponding holomorphic vector fields ``kill them''.  
So in \cite{GMO13} we define the Hirschowitz annihilator $\mathcal A$ to be the connected Lie subgroup of $G$ whose Lie algebra is given by     
\[   
     \mathfrak a \; := \; \{ \ \xi\in\mathfrak g \ |   
    \ \zeta \mapsto \varphi(\exp(\xi\zeta)\cdot x_0) \mbox{ is constant on a neighbourhood of 0 in } \mathbb C \ \forall \varphi\in\mathscr{P}(X) \ \} 
\]  
where $\mathscr P(X)$ denotes the space of continuous plurisubharmonic functions on $X := G/H$.   
In \cite[(3.1)]{GMO13} this subalgebra is denoted by $\widehat{\mathfrak h}$. 
It is defined there more generally for a domain $p:D\to X$ spread over $X$, but here we only need the case $D=X$, $p={\rm id}$. 
The definition in \cite{GMO13} requires that the derivative of $f(t):=\varphi(\exp(t\xi)\cdot x_0)$ with respect to the 
complex variable $t$ vanish near $t=0$ in the sense of distributions, for every continuous psh function $\varphi$. 
Since $\varphi$ is real valued, we have $\partial f/\partial\bar t=\overline{\partial f/\partial t}$. 
Thus $\partial f/\partial t=0$ near $0$ exactly when $df=0$ near $0$, that is, when the continuous function $f$ is constant near $0$, 
and the two definitions give the same subalgebra.
If $\zeta\mapsto\varphi(\exp(\xi\zeta)\cdot x_0)$ is constant near $0$ for every $\varphi\in\mathscr P(X)$, then it is in fact constant on all of $\mathbb C$, 
as one sees by applying the condition to $\varphi\circ\exp(\xi\zeta_0)\in\mathscr P(X)$ for each $\zeta_0\in\mathbb C$.

\medskip 
Clearly $\mathfrak{h} \subset \mathfrak a$. 
Since $\exp(\xi\cdot i\zeta) = \exp((i\xi)\zeta)$, the condition for $i\xi$ is the same as for $\xi$, and 
in fact $\mathfrak a$ is a complex Lie subalgebra of $\mathfrak g$ which properly contains $\mathfrak h$ whenever $X$ is not Stein \cite[Lemma 3.3]{GMO13}. 
Moreover, $H$ normalizes $\mathfrak a$: if $h\in H$ and $\varphi\in\mathscr P(X)$, then $\varphi\circ h\in\mathscr P(X)$ and 
\[
   \varphi\bigl(\exp({\rm Ad}(h)\xi\,\zeta)\cdot x_0\bigr) \; = \; \varphi\bigl(h\exp(\xi\zeta)h^{-1}\cdot x_0\bigr) \; = \; (\varphi\circ h)\bigl(\exp(\xi\zeta)\cdot x_0\bigr),
\]
so that ${\rm Ad}(h)\mathfrak a=\mathfrak a$. 
In particular $H\subset N_G(\mathcal A)$; compare \cite[Proposition 3.6]{GMO13}. 
The same computation, with conjugation by $h$ replaced by a holomorphic automorphism $\alpha$ of $G$ with $\alpha(H)=H$, 
gives $\alpha(\mathcal A)=\mathcal A$. 
In this sense $\mathfrak a$ and $\mathcal A$ are canonically attached to the pair $(G,H)$.

\medskip 
The connected complex subgroup $\mathcal A$ of $G$ need not be closed, but it contains the identity component $H^0$ of $H$. 
Since $H$ normalizes $\mathcal A$, the sets $F_{gx_0}:=g\mathcal A\cdot x_0$, $g\in G$, are well defined. 
The set $F_{gx_0}$ is the orbit through $gx_0$ of the conjugate group $g\mathcal Ag^{-1}$; if $\mathcal A$ is normal in $G$, it is simply the $\mathcal A$-orbit of $gx_0$. 
These sets are the leaves of a complex foliation $\mathscr{F}$, the {\bf Levi foliation} of the manifold $X$.   
By \cite[Thm 3.1]{GMO13} every leaf is a relatively compact immersed complex submanifold of $X$, and we shall use this repeatedly. 
Every $\varphi\in\mathscr P(X)$ is constant on the leaves. 
To see this, recall that $\mathcal A$ is generated by the subgroups $\exp(\mathbb C\xi)$, $\xi\in\mathfrak a$, and that $\varphi\circ a\in\mathscr P(X)$ for $a\in\mathcal A$. 
The definition of $\mathfrak a$ then shows that $\varphi$ is constant on $\mathcal A\cdot x_0$, and applying this to $\varphi\circ g$, $g\in G$, 
gives the same for every leaf. 
In general, complex foliations are rather difficult to understand.      
But here the foliation arises from a subgroup strongly reflecting the underlying geometry of the homogeneous manifold   
and is related to the existence of a plurisubharmonic exhaustion on it.  
This allows a sufficiently good understanding in 
order to analyze the structure from the point of view of its holomorphic reduction.                 

\medskip   
If the leaves of the foliation are closed, then they are compact and $X$ is holomorphically convex \cite[Theorem 4.1]{GMO13}.  
In that case the holomorphic reduction is given by the Remmert reduction \cite{Rem56} and has compact fiber and Stein base.   
The main difficulty occurs when this is not the case.    

\begin{remark} 
In passing we recall that a question of Serre \cite{Ser53} asks whether the total space of every holomorphic fibration with 
Stein fiber and Stein base is itself Stein.  
Counterexamples are known, and some of them are even complex solvmanifolds. 
Such examples are constructed in \cite[\S 4, p. 74]{Loe85}; their holomorphic reduction is a $\mathbb C^* \times \mathbb C^*$-bundle over $\mathbb C^*$, 
although this is not stated explicitly there. 
Lemma \ref{Hirsch} below shows, among other things, that Serre's question has an affirmative answer in the present setting.        
\end{remark}

\begin{lemma} [\cite{Hir75} and Remark 2.5 in \cite{GMO13}]  \label{Hirsch}   
Suppose $Y$ is a pseudoconvex complex homogeneous manifold that admits a holomorphic fibration  
\[  
              Y \; \stackrel{F}{\longrightarrow} \; B  
\]
with both $F$ and $B$ being holomorphically separable manifolds.   
Then $Y$ is Stein.   
\end{lemma}  

\begin{proof}  
Recall that every complex homogeneous manifold is infinitesimally homogeneous.  
Since neither $F$ nor $B$ can contain an inner integral curve, neither can $Y$, 
and so $Y$ is Stein by \cite[Proposition 3.4]{Hir75}.    
\end{proof}  

\subsection{Tits' Normalizer Fibration} \label{titssec}   

Suppose $H$ is a closed complex subgroup of complex dimension $k$ of a complex Lie group $G$ of complex dimension $n$.  
The Lie algebra $\mathfrak h$ of $H$ is a Lie subalgebra of $\mathfrak g$ and can be considered   
as a point in the Grassmannian $Gr(k,n)$ of complex $k$-dimensional vector subspaces of $\mathfrak g$.  
Since ${\rm Ad}(G) \subset GL(\mathfrak g)$, there is a natural action of ${\rm Ad}(G)$ on $Gr(k,n)$ and 
the ${\rm Ad}(G)$-orbit of the point $\mathfrak h$ can be identified with $G/N$, where $N := N_G(H^0)$, 
i.e., the normalizer in $G$ of the identity component $H^0$ of $H$.  
Via the Pl\"ucker embedding $Gr(k,n)\subset\mathbb P(\Lambda^k\mathfrak g)$, on which $G$ acts linearly through $\Lambda^k{\rm Ad}$,   
$G/N$ is realized as an orbit of a linear action of $G$ on a projective space.     
The orbit map $G/N\to\mathbb P(\Lambda^k\mathfrak g)$ is an injective holomorphic immersion. 
If $G$ is not algebraic, however, its image need not be locally closed, and so the manifold topology of $G/N$ may be strictly finer 
than the topology induced from the projective space. 
For example, let $t\in\mathbb C$ act on $\mathbb C^3$ by ${\rm diag}(1,e^t,e^{\alpha t})$, where $\alpha\in\mathbb R\setminus\mathbb Q$, 
let $G:=\mathbb C\ltimes\mathbb C^3$ and $H:=\exp(\mathbb C v)$ with $v=(1,1,1)\in\mathbb C^3$. 
Then the ${\rm Ad}(G)$-orbit of $\mathfrak h$ is $\{[1:e^t:e^{\alpha t}] : t\in\mathbb C\}\cong\mathbb C$, whose closure contains the 
real $3$-dimensional set $\{[1:z:w] : |w|=|z|^\alpha\}$.  
Since $H \subset N_G(H^0)$, we have the Tits normalizer fibration $G/H \to G/N$  \cite{Tits62}.      

\begin{remark}  
Let $I$ be a complex Lie subgroup of the complex Lie group $G$.   
The normalizer $N_G(I^0)$, where $I^0$ denotes the connected component 
of the identity of $I$, is a closed subgroup of $G$, since it is the 
isotropy subgroup of a point in a complex projective space under the adjoint action of the group $G$.  
\end{remark}

\subsection{Fibrations over Projective orbits}  

We use the following notation for the derived series of the Lie group $G$ 
\[  
     G^{(0)} := G, \; G^{(1)} := G' = [G,G], \; \ldots, \; G^{(k)} := [G^{(k-1)},G^{(k-1)}] \mbox{ for all } k \ge1   .   
\]

\begin{lemma}  \label{flag3}  
Let $G$ be a complex Lie group, not necessarily connected, $H$ a closed complex subgroup of $G$, and $I$ a closed complex subgroup of $G$ containing    
$H$ with $G/I$ equivariantly realized as an orbit in the complex projective space $\mathbb P_N$ (as in \S\ref{titssec}).     
Suppose that the fibration $G/H \to G/I$ is a covering.           
If $J$ is a connected, normal, complex Lie subgroup of $G$ whose orbits in $G/H$ are relatively compact, 
then the $J$-orbits in $G/H$ are flag manifolds.   
\end{lemma}  

\begin{proof}  
Let $\rho: G \rightarrow PGL_{N+1}(\mathbb C)$ be the action on $\mathbb P_N$. 
We identify $G/I$ with a $G$-orbit in $\mathbb P_N$. 
The inclusion $G/I\hookrightarrow\mathbb P_N$ is continuous, but we always give $G/I$ its own manifold topology, which may be finer than 
the subspace topology (\S\ref{titssec}). 
Since $G/H \rightarrow G/I$ is continuous, the $J$-orbits in $G/I$ are relatively compact. 
Fix $x\in G/I$ and let $M$ be the closure of $J.x$ in $G/I$. 
Then $M$ is compact, so its image in $\mathbb P_N$ is closed, and hence it is the closure of $J.x$ in $\mathbb P_N$. 

\medskip 
Let $B$ be the Zariski closure of $\rho(J)$ in $PGL_{N+1}(\mathbb C)$; it is a connected linear algebraic group. 
For $g\in G$, conjugation by $\rho(g)$ is an algebraic automorphism that preserves $\rho(J)$, 
because $J$ is a normal subgroup of $G$. 
Hence $\rho(G)$ normalizes $B$ and also its commutator subgroup $B'$. 
By Chevalley's Theorem \cite[Ch. II, \S 14]{Chev51}, one has $B' = [\rho(J),\rho(J)]$ and, in particular 
$B' \subset \rho(J)$. 

\medskip 
We need the following observation. 
Let $E$ be a connected algebraic subgroup of $PGL_{N+1}(\mathbb C)$ that is 
normalized by $\rho(G)$. 
For $y = \rho(g)x$ one has $E.y = \rho(g) (E.x)$ and 
all the $E$-orbits through the points of $G/I$ have the same dimension. 
On the other hand, $E.x$ is Zariski open in its Zariski closure, and the boundary is a union of 
$E$-orbits of strictly smaller dimensions \cite[\S 1.8]{Bor91}. 
The Euclidean closure of a subset of $E.x$ lies in its Zariski closure. 
Hence 
\[ 
   (\ast)\qquad   \mbox{ if } S \subset E.x \mbox{ and the closure of } S \mbox{ in } \mathbb P_N \mbox{ lies in } 
     G/I, \mbox{ then it lies in } E.x 
\] 

\medskip 
Applying $(\ast)$ with $E=B$ and $S=J.x$ gives $M \subset B.x$. 
Since $B'$ is normal in $B$, the group $B_xB'$ is a closed algebraic subgroup and 
the quotient $B/B_xB' \cong (B/B')/(B_xB'/B')$ is a commutative affine algebraic group and, 
as a consequence,  is an affine variety. 
The holomorphic map $J \rightarrow B/B_xB', j \mapsto \rho(j)B_xB'$ factors through $J.x\subset M\subset B.x\cong B/B_x$ 
and the image of the compact set $M$ in the affine variety $B/B_xB'$ is compact.   
Thus its coordinate functions are bounded, and their restrictions to the complex one-parameter subgroups of $J$ 
are bounded entire functions, hence constant.
Consequently the map is constant near $e$, and by the identity principle it is constant on $J$. 
Thus $J.x$ lies in the fiber $B_xB'.x= B'.x$. 
Conversely, $B'.x \subset J.x$ because $B' \subset \rho(J)$.  
Thus $J.x = B'.x$. 

\medskip
Applying $(\ast)$ again with $E=B'$ and $S=J.x$ gives $M \subset B'.x = J.x$. 
Therefore $J.x$ is closed in $\mathbb P_N$ and $J.x = B'.x$ is a compact orbit of a connected linear algebraic group 
whose isotropy is thus parabolic, and $J.x$ is a flag manifold \cite[\S 11.2]{Bor91}. 
Since the bijective $J$-equivariant holomorphic map from $J/J_x$ onto $B'.x$ is biholomorphic, the $J$-orbit through $x$ in $G/I$, 
with its own manifold structure, is this flag manifold. 
Finally, the $J$-orbits in $G/H$ are connected coverings of the $J$-orbits in $G/I$ and thus are flag manifolds, 
because flag manifolds are simply connected. 
We remark that the connectedness of $G$ was not used in this proof. 
\end{proof}

\begin{cor}  \label{flag2}  
Let $G$ be a complex Lie group, not necessarily connected, $H$ a closed complex subgroup of $G$, and $I$ a closed complex subgroup of $G$ 
containing $H$ with $G/I$ equivariantly realized as an orbit in the complex projective space $\mathbb P_N$ (as in \S\ref{titssec})   
and the fibers of the fibration $G/H \to G/I$ holomorphically separable.        
If $J$ is a connected, normal, complex Lie subgroup of $G$ whose orbits in $G/H$ are relatively compact, then the $J$-orbits in $G/H$ are flag manifolds.   
\end{cor}

\begin{proof} 
Let $\pi:G/H\to G/I$ be the fibration, fix $x\in G/H$ and put $y:=\pi(x)$.
The orbit $J.y=\pi(J.x)$ is relatively compact in $G/I$.
Lemma \ref{flag3}, applied with $H$ replaced by $I$ (so that the covering is the identity), shows that $J.y$ is a flag manifold.
Now consider $\pi|_{J.x}:J.x\to J.y$.
Its fiber through $x$ is the orbit $J_y.x$ of the isotropy group $J_y$, a closed complex subgroup of $J$.
This orbit lies in the fiber $\pi^{-1}(y)\cong I/H$, which is closed, so $J_y.x$ is relatively compact in $\pi^{-1}(y)$, 
being contained in the compact set ${\rm cl}(J.x)\cap\pi^{-1}(y)$.
Every holomorphic function on $\pi^{-1}(y)$ is therefore bounded on the orbits of the complex one-parameter subgroups of $J_y^0$ through $x$, 
hence constant on them by Liouville's theorem.
Since $\pi^{-1}(y)$ is holomorphically separable, $J_y^0$ fixes $x$.
Hence the fibers of $\pi|_{J.x}$ are discrete, and $J.x\to J.y$ is a connected covering of a flag manifold.
Flag manifolds are simply connected, so $J.x\cong J.y$ is a flag manifold. 
\end{proof}

\subsection{Normality under closure and complexification}  \label{normsec}
For $J$ a connected Lie subgroup of $G$ the complexification $J^{\mathbb C}$ of $J$ is the connected Lie subgroup of $G$ 
corresponding to the Lie algebra $\mathfrak j + i \mathfrak j$, where $\mathfrak j$ denotes the Lie algebra of $J$.        
We shall need to know how normality behaves under closure and complexification.   

\begin{lemma}  \label{normal}  
Suppose $I$ is a connected normal Lie subgroup of a connected Lie subgroup $J$ of a connected complex Lie group $G$.  
Then 
\begin{enumerate} 
\item $I \lhd {\rm cl}_G(J)$  
\item $I^{\mathbb C} \lhd J^{\mathbb C}$
\end{enumerate}  
In particular, $I^{\mathbb C} \lhd \widetilde{J}$, where $\widetilde{J}$ is the smallest connected closed complex subgroup 
of $G$ that contains $J$.  
\end{lemma}   

\begin{proof} 
(1) If $\lim g_n = g \in {\rm cl}_G(J)$, where $g_n \in J$, then $\rm{Ad}(g_n)\mathfrak{i} = \mathfrak{i}$ for all $n$. 
As a consequence, ${\rm Ad}(g)\mathfrak{i} = \lim {\rm Ad}(g_n)\mathfrak{i} = \mathfrak{i}$, 
since the stabilizer of $\mathfrak{i}$ in $GL(\mathfrak{g})$ is closed. 
As $I$ is connected, $gIg^{-1} = I$. 

\medskip\noindent
(2) This follows from the fact that 
$ [\mathfrak i, \mathfrak j ] \subset \mathfrak i \Longrightarrow [ \mathfrak i + i \mathfrak i , \mathfrak j + i \mathfrak j ] 
\subset \mathfrak i + i \mathfrak i$.  

\medskip\noindent
Finally, $\widetilde{J}$ is obtained from $J$ by alternately taking complexifications and closures, 
and applying (1) and (2) at each stage completes the proof.  
\end{proof}

The group $L$ in Section 3 is not connected, and so we also need an analogue of $\widetilde J$ for groups of the form $L=L^0\Gamma$. 
For a connected Lie subgroup $J$ of $G$ we write $\langle J\rangle$ for the smallest connected closed complex subgroup of $G$ containing $J$. 
This is the identity component of the intersection of all closed complex subgroups of $G$ containing $J$; 
its Lie algebra is the intersection of the corresponding complex subalgebras. 
By construction, every holomorphic automorphism of $G$ that preserves $J$ also preserves $\langle J \rangle$.

\begin{lemma}\label{Ltilde}
Let $\Gamma$ be a discrete subgroup of $G$ and $L$ a closed subgroup with $L=L^0\Gamma$. 
Put $M_0:=\langle L^0\rangle$ and $M_{k+1}:=\langle{\rm cl}_G(M_k\Gamma)^0\rangle$. 
Then each $M_k$ is normalized by $\Gamma$, $M_k\subset M_{k+1}$, and the sequence is eventually constant, say equal to $M_*$. 
The group $\widetilde L:=M_*\Gamma$ is a closed complex subgroup of $G$ with $\widetilde L^0=M_*$ and $L\subset\widetilde L$. 
If $I$ is a connected normal subgroup of $L^0$ normalized by $\Gamma$, then $I^{\mathbb C}\lhd\widetilde L$.
\end{lemma}

\begin{proof}
Conjugation by $\gamma\in\Gamma$ preserves $L^0$, hence $M_0$. 
If $\Gamma$ normalizes $M_k$, then it normalizes the group $M_k\Gamma$, its closure, the identity component of the closure, 
and thus also $M_{k+1}$. 
Clearly $M_k\subset{\rm cl}_G(M_k\Gamma)^0\subset M_{k+1}$. 
Since $\dim M_k \le \dim G$ for all $k$, the increasing sequence of connected groups $M_k$ is eventually constant.  
Writing $M_*$ for the final group, one has   
${\rm cl}_G(M_*\Gamma)^0\subset M_*\subset{\rm cl}_G(M_*\Gamma)$. 
Hence $M_*={\rm cl}_G(M_*\Gamma)^0$ is open in ${\rm cl}_G(M_*\Gamma)$. 
If $x\in{\rm cl}_G(M_*\Gamma)$, the open set $xM_*$ meets $M_*\Gamma$, so $x\in M_*\Gamma M_*=M_*\Gamma$. 
Hence $M_*\Gamma$ is closed and contains $L^0\Gamma=L$. 
For the last statement, $I^{\mathbb C}\lhd M_0$ by Lemma \ref{normal}, and $\Gamma$ normalizes $I^{\mathbb C}$, 
so $I^{\mathbb C}\lhd M_0\Gamma$. 
By (1) of Lemma \ref{normal}, applied with the closed stabilizer of $\mathfrak i+i\mathfrak i$, 
$I^{\mathbb C}$ is normal in ${\rm cl}_G(M_0\Gamma)$, 
and Lemma \ref{normal} gives $I^{\mathbb C}\lhd M_1$.  
Induction gives $I^{\mathbb C}\lhd M_*\Gamma$.
\end{proof}

\subsection{Existence of a fibration by a solvmanifold}   

The following extension of  \cite[Thm 2]{Gi81} yields a fibration with \emph{a complex solvmanifold} as its fiber.  

\begin{proposition}[personal communication from K. Oeljeklaus] \label{Karl}    
Let $X := G/\Gamma$, where $\Gamma$ is a discrete subgroup of a connected, simply connected,   
complex Lie group $G$ with Levi decomposition $G = S \ltimes R$ and $\dim R > 0$.       
Then there is a connected, complex, solvable subgroup $H$ of $G$ normalized by
$\Gamma$, containing $R$, with $H\cdot \Gamma$ a closed subgroup of $G$.   
Moreover, $H$ is closed, and it can be chosen canonically, i.e., such that every holomorphic 
automorphism $\alpha$ of $G$ with $\alpha(\Gamma) = \Gamma$ satisfies $\alpha(H) = H$. 
In particular, one has the non-trivial fibration (unless $S = \{ e \}$)   
\begin{eqnarray} \label{radfibr}  
       G/\Gamma \;  \longrightarrow \;  G/H\cdot \Gamma \; = \; S/S\cap H\cdot \Gamma   
\end{eqnarray}  
that has the connected complex solvmanifold $H/(H \cap \Gamma)$ as its typical fiber.    
\end{proposition}

\begin{proof}  
If the $R$-orbits themselves are closed, set $H := R$.    
If not, then the Zassenhaus Lemma \cite{Aus63} is used in \cite[Theorem 2]{Gi81} to show the existence of    
a connected complex solvable subgroup of $G$ normalized by $\Gamma$ 
and containing the identity component of ${\rm cl}_G(R\cdot\Gamma)$.   
Let $H_1$ be the smallest such subgroup, i.e., the connected subgroup whose Lie algebra is the intersection 
of the Lie algebras of all of them; it is again connected, complex, solvable, normalized by $\Gamma$, 
and contains ${\rm cl}_G(R\cdot\Gamma)^{0}$. 
If $H_1\cdot\Gamma$ is closed in $G$, one has the desired result with $H := H_1$.    

\medskip 
Otherwise, let $N_1 := N_G(H_1)$. 
Since $\Gamma$ normalizes the identity component $N_1^0$ of $N_1$, it also normalizes its radical $R_1$.    
Now $H_1 \subset R_1$, because $H_1$ is solvable and normal in $N_1$.    
If $R_1\cdot\Gamma$ is closed in $G$, we take $H := R_1$. 
Otherwise, note that $\Gamma\subset N_1$, so ${\rm cl}_G(R_1\cdot\Gamma) \subset N_1$, and put $\Gamma_1:=\Gamma\cap N_1^0$. 
Since $N_1^0$ is open in $N_1$ and $R_1\Gamma\cap N_1^0=R_1\Gamma_1$, one has 
${\rm cl}_G(R_1\Gamma)\cap N_1^0={\rm cl}_{N_1^0}(R_1\Gamma_1)$, and so the identity components of ${\rm cl}_G(R_1\Gamma)$ and 
${\rm cl}_{N_1^0}(R_1\Gamma_1)$ coincide. 
The Zassenhaus Lemma \cite{Aus63}, applied to the discrete subgroup $\Gamma_1$ of the connected group $N_1^0$ with radical $R_1$, 
shows that this identity component is solvable; it is normalized by $\Gamma$, because $\Gamma$ normalizes $R_1$.  
Let $H_2$ be the smallest connected closed complex subgroup of $G$ that contains this identity component.   
Then $H_2$ is solvable and normalized by $\Gamma$, and $H_1\subset R_1\subset{\rm cl}_G(R_1\Gamma)^0\subset H_2$. 
If $\dim H_2=\dim H_1$, then all these connected groups coincide, so $R_1={\rm cl}_G(R_1\Gamma)^0$ is open in ${\rm cl}_G(R_1\Gamma)$, 
and, as in the proof of Lemma \ref{Ltilde}, $R_1\Gamma$ is closed, contrary to our assumption. 
Hence $\dim H_2>\dim H_1$, and we repeat the construction with $H_2$ in place of $H_1$. 
A finite number of steps yields the desired connected complex solvable group $H$.     
Each operation used in the construction (taking radicals, closures, identity components and normalizers, and passing to the 
smallest subgroup with a given property) commutes with automorphisms. 
Hence a holomorphic automorphism $\alpha$ of $G$ with $\alpha(\Gamma)=\Gamma$ preserves each of the groups that occur, 
and in particular $H$. 
Finally, $H\Gamma$ is closed and $H\Gamma/H$ is countable, so by Baire's theorem  the connected Lie subgroup $H$ is open 
in $H\Gamma$. 
Hence $H = (H\Gamma)^0$ is closed. 
\end{proof}   

\subsection{Closedness of complexified orbits in Stein solvmanifolds}

\begin{lemma}\label{open}
Let $M$ be a topological group, $U\subset M$ an open subgroup and $P\subset M$ a subgroup.
If $P\cap U$ is closed in $U$, then $P$ is closed in $M$.
\end{lemma}

\begin{proof}
An open subgroup is also closed, so $P\cap U$ is closed in $M$.
Let $x\in{\rm cl}_M(P)$ and $p_k\in P$ with $p_k\to x$.
Then $p_k\in xU$ for $k\ge k_0$, so $p_{k_0}^{-1}p_k\in P\cap U$ converges to $p_{k_0}^{-1}x$, which therefore lies in $P\cap U$.
Hence $x\in P$.
\end{proof}

In a simply connected solvable Lie group every connected Lie subgroup is closed and simply connected \cite{Che41}, 
and the quotient by a connected normal subgroup is again simply connected and solvable. 

\begin{proposition}\label{abel}
Let $H$ be a simply connected solvable complex Lie group, $\Delta\subset H$ a discrete subgroup with $H/\Delta$ Stein.
Let $Z\subset H$ be a connected abelian real subgroup with $\mathfrak z\cap i\mathfrak z=0$.
Suppose $Z^{\mathbb C}$ is normal in $H$, $\Delta$ normalizes $Z$, and $Z\cap\Delta$ is a lattice in $Z$.
Then:
\begin{enumerate}[\rm(i)]
\item $Z^{\mathbb C}\cap\Delta=Z\cap\Delta$;
\item $Z^{\mathbb C}\Delta$ is closed in $H$, 
  so $Z^{\mathbb C}\Delta/\Delta\cong Z^{\mathbb C}/(Z\cap\Delta)\cong(\mathbb C^*)^k$, $k=\dim_{\mathbb R}Z$;
\item $H/Z^{\mathbb C}\Delta$ is Stein.
\end{enumerate}
\end{proposition}

\begin{proof}
Identify $Z^{\mathbb C}$ with $\mathfrak z^{\mathbb C}=\mathfrak z\oplus i\mathfrak z\cong\mathbb C^k$ via $\exp$, 
and put $\Lambda_0 := Z\cap\Delta$.  
Then $\Lambda_0$ is a lattice in $\mathfrak z\cong\mathbb R^k$.

\medskip 
\emph{Step A: Reduction to the case where $Z^{\mathbb C}$ is central.}  
Let $\rho:H\to GL(\mathfrak z^{\mathbb C})$ be the adjoint action on the abelian ideal $\mathfrak z^{\mathbb C}$, 
and $H_1:=\ker\rho$, the centralizer of $Z^{\mathbb C}$.
It is a closed complex subgroup containing $Z^{\mathbb C}$.
Each $\delta\in\Delta$ normalizes $Z$ and $\Delta$, so $\rho(\delta)$ preserves $\mathfrak z$ and $\Lambda_0$.
Hence $\rho(\Delta)$ lies in the stabilizer of $\Lambda_0$, which is a copy of $GL_k(\mathbb Z)$ and is discrete in $GL(\mathfrak z^{\mathbb C})$.
Let $p:H\to H/H_1$.
If $p(\delta_m)\to e$ then $\rho(\delta_m)\to 1$, so $\rho(\delta_m)=1$ for large $m$.
Thus $p(\Delta)$ is discrete, $H_1\Delta=p^{-1}(p(\Delta))$ is closed, and $H_1$ is open in $H_1\Delta$.
Hence $H_1^0$ is open in $H_1^0\Delta$, and $H_1^0\Delta$ is closed, being an open subgroup of the closed group $H_1\Delta$.

\medskip 
Put $\Delta_1:=\Delta\cap H_1^0\supset\Lambda_0$.
Then $Y_1:=H_1^0/\Delta_1\cong H_1^0\Delta/\Delta$ is a closed submanifold of $H/\Delta$, and hence it is Stein.
Also $Z^{\mathbb C}\Delta\cap H_1^0=Z^{\mathbb C}\Delta_1$ and $Z^{\mathbb C}\cap\Delta=Z^{\mathbb C}\cap\Delta_1$.
By Lemma \ref{open}, applied in $M=H_1^0\Delta$ with the open subgroup $U=H_1^0$, 
it suffices to prove (i) and (ii) for $(H_1^0,\Delta_1,Z)$, where $Z^{\mathbb C}$ is central.

\medskip 
\emph{Step B: The action is free, and (i) holds.}  
Since $Z^{\mathbb C}$ is central in $H_1^0$, left multiplication gives a holomorphic action of $Z^{\mathbb C}$ on $Y_1$.
The stabilizer of every point is $\Lambda:=Z^{\mathbb C}\cap\Delta_1$.
The orbit map $Z^{\mathbb C}/\Lambda\to Y_1$ is an injective holomorphic immersion into a Stein manifold, 
so $Z^{\mathbb C}/\Lambda$ is holomorphically separable.
Let $V$ be the real span of $\Lambda$.
Since $\Lambda$ is discrete and spans $V$, it is a lattice in $V$, and $V/\Lambda$ is compact.
For $e\in V\cap iV$ and any holomorphic $f$ on $\mathbb C^k/\Lambda$, the entire function $\zeta\mapsto f(p+\zeta e)$ 
takes values in the image of the compact set $p+V/\Lambda$, so it is constant.
Hence $V\cap iV=0$.
But $V\supset\mathfrak z$, because $\Lambda\supset\Lambda_0$, and a totally real subspace of $\mathfrak z\oplus i\mathfrak z$ 
has real dimension at most $k$.
Therefore $V=\mathfrak z$ and $\Lambda=\Lambda_0$, which proves (i).
In particular the torus $T^{\mathbb C}:=Z^{\mathbb C}/\Lambda_0\cong(\mathbb C^*)^k$ acts freely on $Y_1$, 
and its compact real form $T:=Z/\Lambda_0$ acts as well.

\medskip
\emph{Step C: The orbits are closed, and (ii) holds.}  
Let $\varrho$ be a smooth strictly psh exhaustion of $Y_1$, and $\varrho_T$ its average over $T$.
It is a smooth, strictly psh, $T$-invariant exhaustion.
For $y\in Y_1$ and $s\in\mathfrak z$ put $v_y(s):=\varrho_T(\exp(is)\cdot y)$.
For a unit vector $\theta\in\mathfrak z$, the function $\zeta\mapsto\varrho_T(\exp(\zeta\theta)\cdot y)$ is subharmonic on $\mathbb C$ 
and depends only on ${\rm Im}\,\zeta$, by $T$-invariance.
Hence $g(t):=v_y(t\theta)$ satisfies
\[
      g''(t) = 4\,\mathcal L_{\varrho_T}(\theta_Y)\bigl(\exp(it\theta)\cdot y\bigr),
\]
where $\theta_Y$ is the holomorphic vector field of $\theta$ on $Y_1$ and $\mathcal L$ is the Levi form.
The field $\theta_Y$ has no zeros, because the action is free.
On the compact set $K_c:=\{\varrho_T\le c\}$ the continuous function $(\theta,y)\mapsto\mathcal L_{\varrho_T}(\theta_Y)(y)$, $|\theta|=1$, 
has a positive minimum $\varepsilon_c$.

\medskip  
We claim that $v_y$ is proper.
Otherwise there are $s_m\to\infty$ with $v_y(s_m)\le c$.
Passing to a subsequence, $s_m/|s_m|\to\theta$.
By convexity, $v_y\le\max(v_y(0),c)=:c'$ on each segment $[0,s_m]$, and by continuity $g(t)=v_y(t\theta)\le c'$ for all $t\ge0$.
So $\exp(it\theta)\cdot y\in K_{c'}$ for $t\ge0$, and $g''\ge4\varepsilon_{c'}>0$ on $[0,\infty)$.
Then $g(t)\to\infty$, a contradiction.

\medskip 
Now let $z_m\cdot y\to y'$ with $z_m\in T^{\mathbb C}$.
Write $z_m=\exp(is_m)\tau_m$ with $s_m\in\mathfrak z$ and $\tau_m\in T$.
Then $v_y(s_m)=\varrho_T(z_m\cdot y)$ is bounded, so $(s_m)$ is bounded, and after passing to a subsequence $z_m\to z$.
Thus $y'=z\cdot y$, and every $T^{\mathbb C}$-orbit is closed.
For the orbit of the base point this says that $Z^{\mathbb C}\Delta_1$ is closed in $H_1^0$, which proves (ii).

\medskip 
\emph{Step D: The base is Stein.} 
Here we work on $H$ itself, where $Z^{\mathbb C}$ need not be central.
Let $\tilde\varrho$ be the lift to $H$ of a smooth strictly psh exhaustion of $H/\Delta$, and put
\[
\varrho'(h):=\int_{Z/\Lambda_0}\tilde\varrho(hz)\,dz .
\]
This is well defined on $H/\Delta$, since for $\delta\in\Delta$ we have $\tilde\varrho(h\delta z)=\tilde\varrho(h\,\delta z\delta^{-1})$.
Moreover $z\mapsto\delta z\delta^{-1}$ is an automorphism of $Z$ preserving $\Lambda_0$, so it has determinant $\pm1$ and preserves Haar measure.
The function $\varrho'$ is smooth, strictly psh (an average of right translates of $\tilde\varrho$), right $Z$-invariant, and an exhaustion. 
To see the last property, let $C_Z\subset Z$ be a compact set with $C_Z\Lambda_0=Z$, and $C\subset H$ a compact set mapping onto $\{\varrho\le c\}$. 
If $\varrho'(h\Delta)\le c$, then $\varrho(hz\Delta)\le c$ for some $z\in Z$, so $hz\in C\Delta$ and we may take $z\in C_Z$.  
Hence $h\in C\Delta C_Z^{-1}\subset CZ\Delta$, because $\Delta$ normalizes $Z$, and $CZ\Delta=CC_Z\Lambda_0\Delta=CC_Z\Delta$. 
So $h\Delta$ lies in the compact set $CC_Z\Delta/\Delta$.

\medskip 
Let $q:H/\Delta\to B:=H/Z^{\mathbb C}\Delta$, choose a local holomorphic section $\sigma:U\to H$ over a ball $U\subset B$, and 
put $F(u,w):=\varrho'(\sigma(u)\exp(w))$ for $w=a+ib\in\mathfrak z^{\mathbb C}$.
Then $F$ is smooth, strictly psh, and independent of $a$.
By (i) the fiber $q^{-1}(u)$ is $\mathfrak z^{\mathbb C}/\Lambda_0$ in the coordinate $w$, and by (ii) it is closed.
Consequently $b\mapsto F(u,ib)$ is proper and strictly convex. 
Its Hessian is $4$ times the Levi form in the fiber directions.
It therefore has a unique, non-degenerate minimum $b(u)$, which depends smoothly on $u$ by the implicit function theorem.
Put $m(u):=F(u,ib(u))=\min q^{-1}(u)\,\varrho'$, a smooth function on $B$, independent of the section.

\medskip 
Fix $u_0$, a small ball $U_0\ni u_0$, and a small convex ball $\omega$ around $b(u_0)$ with $b(U_0)\subset\omega$.
On the relatively compact set $U_0\times(\mathfrak z+i\omega)$ (modulo $\Lambda_0$), 
the function $F-\varepsilon|u|^2$ is psh for small $\varepsilon>0$.
By Kiselman's minimum principle \cite[Theorem 2.2]{Kis78}, $\inf_{b\in\omega}(F(u,ib)-\varepsilon|u|^2)$ is psh on $U_0$.
By convexity this infimum is $m(u)-\varepsilon|u|^2$.
Hence $m$ is strictly psh.
It is an exhaustion, since $\{m\le c\}=q(\{\varrho'\le c\})$.
Therefore $B$ is Stein, see \cite{Gra58} or \cite{Nar62}, which proves (iii).
\end{proof}

\begin{definition}  
A $\mathbb C^*$ power tower of length {\it zero} is a point, and one of length {\it one} is simply the manifold $(\mathbb C^*)^p$ for some positive integer $p$.   
For any integer $n>1$ a $\mathbb C^*$ power tower of length $n$ is the total space of a holomorphic fiber bundle with fiber $(\mathbb C^*)^k$, 
for some positive integer $k$, over a $\mathbb C^*$ power tower of length $n-1$.      
(The bundles that occur below need not be principal, since $\Delta$ may act non-trivially on the lattice $\Lambda_0$ by conjugation.)
\end{definition}   

\begin{thm}\label{solv}
Let $H$ be a simply connected solvable complex Lie group, $\Delta\subset H$ a discrete subgroup with $H/\Delta$ Stein.
Let $W\subset H$ be a connected closed real subgroup with $\mathfrak w\cap i\mathfrak w=0$.
Suppose $W^{\mathbb C}$ is normal in $H$, $\Delta$ normalizes $W$, and $W\cap\Delta$ is cocompact in $W$.
Then $W^{\mathbb C}\Delta$ is closed in $H$, and $H/W^{\mathbb C}\Delta$ is Stein.
Moreover $W^{\mathbb C}\Delta/\Delta$ is a $\mathbb C^*$ power tower, each of whose stages is Stein.
\end{thm}

\begin{proof}
The proof is by induction on $\dim W$.  
The case $W=\{e\}$ is trivial.
Let $N_W$ be the nilradical of $W$ and $Z:=Z(N_W)$, a non-trivial connected abelian group, characteristic in $W$, and hence normalized by $\Delta$.
Since $W\cap\Delta$ is a lattice in $W$, the group $N_W\cap\Delta$ is a lattice in $N_W$ \cite[Cor.~3.5]{Rag72}, and 
$Z\cap\Delta$ is a lattice in $Z$ \cite[Prop.~2.17]{Rag72}. 
The nilradical of $\mathfrak w^{\mathbb C}$ is $\mathfrak n_W^{\mathbb C}$, and its center is $\mathfrak z^{\mathbb C}$.
Hence $Z^{\mathbb C}$ is characteristic in $W^{\mathbb C}$, hence normal in $H$.
Also $\mathfrak z\cap i\mathfrak z=0$.
By Proposition \ref{abel}, $Z^{\mathbb C}\Delta$ is closed and $H/Z^{\mathbb C}\Delta$ is Stein.

\medskip 
Put $\bar H:=H/Z^{\mathbb C}$, and let $\bar\Delta$ and $\bar W$ be the images of $\Delta$ and $W$.
Because $Z^{\mathbb C}\Delta$ is closed in $H$, $\bar\Delta$ is a countable closed subgroup of $\bar H$ and hence discrete, 
and $\bar H/\bar\Delta=H/Z^{\mathbb C}\Delta$ is Stein.
Also $\bar W^{\mathbb C}=W^{\mathbb C}/Z^{\mathbb C}$ is normal in $\bar H$, and $\bar\Delta$ normalizes $\bar W$.
The image of $W\cap\Delta$ is cocompact in $\bar W$, so the discrete group $\bar W\cap\bar\Delta$ is cocompact in $\bar W$.

\medskip 
Next, $\bar W$ is totally real.
To see this, note first that $\mathfrak w\cap\mathfrak z^{\mathbb C}=\mathfrak z$, because $\mathfrak w\cap i\mathfrak w=0$.
Now let $x\in\mathfrak w$ with $ix\in\mathfrak w+\mathfrak z^{\mathbb C}$, say $ix=y+a+ib$ with $y\in\mathfrak w$ and $a,b\in\mathfrak z$.
Then $i(x-b)=y+a\in\mathfrak w$ and $x-b\in\mathfrak w$, so $x=b\in\mathfrak z$.

\medskip 
By induction, $\bar W^{\mathbb C}\bar\Delta$ is closed and $\bar H/\bar W^{\mathbb C}\bar\Delta$ is Stein.
Its preimage in $H$ is $W^{\mathbb C}\Delta$, and $\bar H/\bar W^{\mathbb C}\bar\Delta=H/W^{\mathbb C}\Delta$.
Each stage is a closed submanifold of the Stein manifold $H/Z_i^{\mathbb C}\Delta$ and 
the statement about the tower follows.   
\end{proof}

\subsection{Invariance and the semisimple case}

As in \S\ref{normsec}, for a connected real Lie subgroup $S$ of a complex Lie group we denote by $S^{\mathbb C}$ 
the connected complex subgroup with Lie algebra $\mathfrak s+i\mathfrak s$.
If $S$ is semisimple, we write $\mathfrak s=\mathfrak s_{nc}\oplus\mathfrak s_c$, where $\mathfrak s_{nc}$ (resp. $\mathfrak s_c$) 
is the sum of the non-compact (resp. compact) simple ideals of $\mathfrak s$, and $S_{nc}$ (resp. $S_c$) denotes the connected Lie subgroup 
of $S$ with Lie algebra $\mathfrak s_{nc}$ (resp. $\mathfrak s_c$).

\begin{proposition}\label{inv}
Let $M$ be a connected complex Lie group and $S\subset M$ a connected semisimple real Lie subgroup.
Let $h: M \to \mathbb R$ be continuous and plurisubharmonic with $h(ms) = h(m)$ for all $m\in M$ and $s\in S$.
Then $h(mc)=h(m)$ for all $c\in S_{nc}^{\mathbb C}$.
In particular, if $S$ has no compact simple factors, $h$ is right invariant under $S^{\mathbb C}$.
\end{proposition}

\begin{proof}
Let $\mathfrak s_j$ denote the non-compact simple ideals of $\mathfrak s$.
For each $j$ we choose a Cartan decomposition $\mathfrak s_j=\mathfrak k_j\oplus\mathfrak p_j$ and a maximal abelian subspace 
$\mathfrak a_j\subset\mathfrak p_j$, which is non-zero since $\mathfrak s_j$ is not compact.
Since the restricted roots are real on $\mathfrak a_j$, we may choose $0\neq X_j\in\mathfrak a_j$ on which every restricted root 
takes a rational value, and we put $X:=\sum_jX_j$; then ${\rm ad}\,X$ is diagonalizable on $\mathfrak s+i\mathfrak s$ with rational eigenvalues.

\medskip
We claim that $t\mapsto\exp(itX)$ is periodic.
The algebra $\mathfrak s+i\mathfrak s$ is a quotient of $\mathfrak s\otimes\mathbb C$, hence semisimple, and so $S^{\mathbb C}$ is a 
connected complex semisimple Lie group and has finite center.
Now ${\rm Ad}(\exp(itX))|_{\mathfrak s+i\mathfrak s}=e^{it\,{\rm ad}X}$ is periodic in $t$, since the eigenvalues of ${\rm ad}\,X$ are rational, 
and so $\exp(iT_0X)$ is central in $S^{\mathbb C}$ for some $T_0>0$.
As the center is finite, the claim follows.

\medskip
Now fix $m\in M$ and consider the subharmonic function $u(s+it):=h(m\exp((s+it)X))$ on $\mathbb C$.
Since $\exp((s+it)X)=\exp(itX)\exp(sX)$ with $\exp(sX)\in S$, one has $u(s+it)=h(m\exp(itX))$, i.e., $u$ depends only on $t$ 
and is therefore a convex function of $t$.
Being periodic and convex, it is constant, and so $h(m\exp(itX))=h(m)$ for all $m\in M$ and $t\in\mathbb R$.

\medskip
It remains to pass from $\exp(i\mathbb R X)$ to all of $S_{nc}^{\mathbb C}$.
Since $h$ is continuous, the set of $c\in M$ with $h(mc)=h(m)$ for all $m\in M$ is a closed subgroup of $M$; let $\mathfrak t$ denote its Lie algebra.
Then $\mathfrak s\subset\mathfrak t$ and $iX\in\mathfrak t$, and $\mathfrak i:=\{Y\in\mathfrak s: iY\in\mathfrak t\}$ is an ideal of $\mathfrak s$, 
because $i[Z,Y]=[Z,iY]\in\mathfrak t$ for $Z\in\mathfrak s$ and $Y\in\mathfrak i$.
As $\mathfrak i$ contains $X$ and every $X_j$ is non-zero, it follows that $\mathfrak i\supset\mathfrak s_{nc}$, i.e., 
$\mathfrak t\supset\mathfrak s_{nc}+i\mathfrak s_{nc}$. 
\end{proof}

\begin{lemma}\label{ssclosed}
Let $M$ be a Lie group and $S\subset M$ an analytic subgroup whose Lie algebra $\mathfrak s$ is semisimple and whose center $Z(S)$ is finite.
Then $S$ is closed in $M$.
\end{lemma}

\begin{proof}
Let $\Sigma:={\rm cl}_M(S)$, a closed connected subgroup of $M$ with Lie algebra $\mathfrak g_\Sigma\supset\mathfrak s$.
By (1) of Lemma \ref{normal}, whose proof works for any Lie group, $S$ is normal in $\Sigma$, so that $\mathfrak s$ is an ideal of $\mathfrak g_\Sigma$.
We claim that $\mathfrak g_\Sigma=\mathfrak s\oplus\mathfrak c$, where $\mathfrak c:=\{X\in\mathfrak g_\Sigma:[X,\mathfrak s]=0\}$ is the centralizer of $\mathfrak s$, 
which is also an ideal.
Indeed, $\mathfrak s\cap\mathfrak c$ is the center of $\mathfrak s$ and so is trivial, and for $X\in\mathfrak g_\Sigma$ the derivation ${\rm ad}\,X|_{\mathfrak s}$ 
of $\mathfrak s$ is inner, i.e., ${\rm ad}\,X|_{\mathfrak s}={\rm ad}\,Y|_{\mathfrak s}$ for some $Y\in\mathfrak s$, so that $X-Y\in\mathfrak c$.

\medskip
Now let $C$ be the analytic subgroup of $\Sigma$ with Lie algebra $\mathfrak c$.
Since $S$ and $C$ are connected and $[\mathfrak s,\mathfrak c]=0$, they commute, and so $q:S\times C\to\Sigma$, $(s,c)\mapsto sc$, 
is a homomorphism of Lie groups (each factor carrying its own Lie group topology).
Its differential is the isomorphism $\mathfrak s\oplus\mathfrak c\to\mathfrak g_\Sigma$, so its image is an open subgroup of the connected group $\Sigma$, 
i.e., $q$ is surjective.
Its kernel $F:=\{(z,z^{-1}):z\in S\cap C\}$ is finite, because $S\cap C\subset Z(S)$, and since a surjective homomorphism between 
Lie groups with countably many components is open, $q$ induces an isomorphism $(S\times C)/F\cong\Sigma$.
Finally, $(S\times\{e\})F$ is a finite union of cosets of the closed subgroup $S\times\{e\}$, and so it is closed and saturated 
with respect to $F$; hence its image $q(S\times\{e\})=S$ is closed in $\Sigma$, and so in $M$. 
\end{proof} 

\begin{lemma}\label{sshull}
Let $M$ be a complex Lie group and $P\subset M$ a closed subgroup whose identity component $S:=P^0$ is semisimple.
Then $Q:=S^{\mathbb C}P$ is a closed subgroup of $M$ with $Q^0=S^{\mathbb C}$.
\end{lemma}

\begin{proof}
A connected complex semisimple Lie group has finite center (it is a quotient of the simply connected one, whose center is finite), 
so $S^{\mathbb C}$ is closed by Lemma \ref{ssclosed}.
Since $P$ normalizes $S$, hence also $\mathfrak s+i\mathfrak s$ and $S^{\mathbb C}$, it follows that $Q$ is a group.

\medskip 
Let $\rho:P\to{\rm Aut}(\mathfrak s)$, $p\mapsto{\rm Ad}(p)|_{\mathfrak s}$.
As ${\rm Aut}(\mathfrak s)/{\rm Int}(\mathfrak s)$ is finite and ${\rm Int}(\mathfrak s)=\rho(S)$, the subgroup $P_1:=\rho^{-1}(\rho(S))$ 
has finite index in $P$.
One has $P_1=S\cdot(P\cap Z)$, where $Z$ denotes the centralizer of $S^{\mathbb C}$ in $M$, and $(P\cap Z)^0\subset S\cap Z\subset Z(S)$, 
which is discrete, so that $P\cap Z$ is discrete.

\medskip
Now $N:=N_M(S^{\mathbb C})$ is a closed subgroup of $M$ whose Lie algebra is $(\mathfrak s+i\mathfrak s)\oplus\mathfrak z_{\mathfrak m}(\mathfrak s+i\mathfrak s)$, 
since the centralizer of a semisimple ideal is a complementary ideal.
Thus $S^{\mathbb C}Z$ contains $N^0$, and so it is open, and hence closed, in $N$.
The homomorphism $Z\to S^{\mathbb C}Z/S^{\mathbb C}$ is surjective with finite kernel $Z\cap S^{\mathbb C}\subset Z(S^{\mathbb C})$, 
and so it induces an isomorphism of Lie groups $Z/(Z\cap S^{\mathbb C})\cong S^{\mathbb C}Z/S^{\mathbb C}$.
Consequently, the image of the discrete group $P\cap Z$ is discrete, i.e., $S^{\mathbb C}(P\cap Z)=S^{\mathbb C}P_1$ is closed.
Since $Q$ is a finite union of cosets of $S^{\mathbb C}P_1$, it is also closed, and $Q/S^{\mathbb C}$ is discrete. 
\end{proof}

\begin{lemma}\label{proper}
Let $T$ be a connected complex semisimple Lie group and $S\subset T$ a connected semisimple real Lie subgroup with $\mathfrak s+i\mathfrak s=\mathfrak t$.
Put $K:=S_c$ and $T_{nc}:=S_{nc}^{\mathbb C}$, and let $p:T\to\bar T:=T/T_{nc}$ be the quotient map.
\begin{enumerate}[\rm(1)]
\item $K$ is compact, $T_{nc}$ is normal in $T$, $T=T_{nc}K^{\mathbb C}$, the subgroups $T_{nc}$ and $K^{\mathbb C}$ commute, 
and $p|_{K^{\mathbb C}}:K^{\mathbb C}\to\bar T$ is a surjective homomorphism with finite kernel; in particular it is proper.
\item If $t\in T$ normalizes $S$, then $p(t)\in p(K)$.
\end{enumerate}
\end{lemma}

\begin{proof}
(1) The ideals $\mathfrak s_{nc}$ and $\mathfrak s_c$ of $\mathfrak s$ commute, so $\mathfrak t_{nc}:=\mathfrak s_{nc}+i\mathfrak s_{nc}$ and 
$\mathfrak k^{\mathbb C}:=\mathfrak s_c+i\mathfrak s_c$ are commuting complex ideals of $\mathfrak t$ with sum $\mathfrak t$. 
Their intersection lies in $\mathfrak k^{\mathbb C}$ and commutes with it, so it lies in the center of the semisimple algebra $\mathfrak k^{\mathbb C}$, 
which is trivial; thus $\mathfrak t=\mathfrak t_{nc}\oplus\mathfrak k^{\mathbb C}$. 
The group $K$ is compact, since its Lie algebra $\mathfrak s_c$ is compact and semisimple, and the connected groups $T_{nc}$ and $K^{\mathbb C}$ commute. 
As $T_{nc}$ is normal and $T_{nc}K^{\mathbb C}$ is a subgroup with Lie algebra $\mathfrak t$, 
we have $T=T_{nc}K^{\mathbb C}$, and so $p|_{K^{\mathbb C}}$ is surjective. 
Its kernel $K^{\mathbb C}\cap T_{nc}$ is a discrete normal subgroup of the connected group $K^{\mathbb C}$ and is therefore central; 
since a connected complex semisimple Lie group has finite center, the kernel is finite, and $p|_{K^{\mathbb C}}$ is a finite covering map, in particular proper.

\medskip\noindent
(2) Since ${\rm Ad}(t)$ preserves $\mathfrak s$, it permutes the simple ideals of $\mathfrak s$ and preserves their compactness, and so it preserves $\mathfrak s_c$. 
Thus $t$ normalizes $K$, and $p(t)$ normalizes $U:=p(K)$.
Moreover, $\mathfrak s_c\cap i\mathfrak s_c=0$, since it is a complex ideal of the compact semisimple Lie algebra $\mathfrak s_c$ and the Killing form 
of a non-zero complex semisimple Lie algebra, regarded as a real Lie algebra, is indefinite. 
So, by (1), $\bar T$ is a connected complex semisimple Lie group with Lie algebra $\mathfrak u\oplus i\mathfrak u$, $\mathfrak u$ the Lie algebra of $U$, 
i.e., $U$ is a compact real form of $\bar T$. 
We claim that $U$ is its own normalizer in $\bar T$. 
Regard $\bar T$ as a real connected semisimple Lie group with Cartan decomposition $\mathfrak u\oplus i\mathfrak u$; its analytic subgroup with 
Lie algebra $\mathfrak u$ is $U$. 
By \cite[Ch.~VI, Theorem 1.1(iii)]{Hel78}, every $g\in\bar T$ can be written as $g=\exp(iX)u$ with $X\in\mathfrak u$ and $u\in U$. 
If $g$ normalizes $U$, then so does $\exp(iX)=gu^{-1}$, and $A:={\rm Ad}(\exp(iX))=e^{i\,{\rm ad}X}$ preserves $\mathfrak u$, i.e., 
it commutes with the conjugation $\tau$ of $\bar{\mathfrak t}$ with respect to $\mathfrak u$. 
Since ${\rm ad}X$ commutes with $\tau$ and $\tau$ is antilinear, $\tau A\tau=e^{-i\,{\rm ad}X}=A^{-1}$, and so $e^{2i\,{\rm ad}X}=1$. 
With respect to the Hermitian inner product $-B(x,\tau y)$, $B$ the Killing form, ${\rm ad}X$ is skew-Hermitian, 
so $e^{2i\,{\rm ad}X}$ is positive definite Hermitian, and $e^{2i\,{\rm ad}X}=1$ forces ${\rm ad}X=0$, i.e., $X=0$. 
Thus $g=u\in U$, and in particular $p(t)\in U=p(K)$.
\end{proof}

In the proof of the Main Theorem we will need the following minimum principle in order to handle 
the compact simple factors of semisimple groups. 
 
\medskip
\noindent\textbf{(MP)}
\emph{Let $K$ be a connected compact Lie group, $K^{\mathbb C}$ a connected complex Lie group with $K$ as a closed real form,
$U\subset\mathbb C^m$ a domain, and $\psi$ a continuous plurisubharmonic function on $U\times K^{\mathbb C}$ with $\psi(u,gk)=\psi(u,g)$ for $k\in K$,
such that for every compact set $C\subset U$ and every $c\in\mathbb R$ the set $\{(u,g)\in C\times K^{\mathbb C}:\psi(u,g)\le c\}$ is compact.
Then $m(u):=\inf_{g\in K^{\mathbb C}}\psi(u,g)$ is plurisubharmonic on $U$.}

\medskip
If $K$ is a torus, this is Kiselman's minimum principle \cite[Theorem 2.2]{Kis78}, which does not even need the exhaustion property.
In general (MP) follows from a theorem of Loeb \cite[Th\'eor\`eme 3]{Loe85}, which is proved in \S 5 of that paper
for pseudoconvex pairs $(G_{\mathbb C}, G_{\mathbb R})$ such that $G_{\mathbb R}$ contains a uniform discrete subgroup.
Here $G_{\mathbb C} = K^{\mathbb C}$ and $G_{\mathbb R} = K$, and since $K$ is compact the trivial group $\{e\}$ is a uniform discrete subgroup.
Loeb's theorem asserts that if $u$ is a $\mathcal C^\infty$ plurisubharmonic function on a tube $X \times K^{\mathbb C}$ that is right $K$-invariant,
strictly plurisubharmonic on each fiber $\{x\}\times K^{\mathbb C}$, and induces an exhaustion of each fiber $\{x\}\times K^{\mathbb C}/K$,
then $x \mapsto \inf_b u(x,b)$ is $\mathcal C^\infty$ and plurisubharmonic on $X$.
In order to deduce (MP), which only assumes continuity, we regularize as follows.

\begin{proof}[Proof of (MP)]
Since plurisubharmonicity is a local property, it suffices to prove that $m$ is plurisubharmonic on every ball $U_0$ with 
$\overline{U_0} \subset U$.   
Choose a ball $U_1$ with $\overline{U_0} \subset U_1$ and $\overline{U_1} \subset U$.   
First, $\psi$ is bounded from below on $\overline{U_1} \times K^{\mathbb C}$, and, as in Lemma \ref{fiberops}, 
the infimum defining $m$ is attained and $m$ is continuous.    
Let $\rho \ge 0$ be a smooth strictly plurisubharmonic exhaustion of the Stein manifold $K^{\mathbb C}$; replacing $\rho$ by its
average $\int_K \rho(gk)\,dk$ over right translations by $K$ we may assume that $\rho$ is right $K$-invariant.
(The existence of such a $\rho$ is also what makes $(K^{\mathbb C},K)$ a pseudoconvex pair in the sense of \cite{Loe85}.)

\medskip
Now let $\chi_\delta \ge 0$, $\delta > 0$, be smooth functions with integral one on $\mathbb C^m \times K^{\mathbb C}$ 
(with respect to the product of
Lebesgue measure and a right Haar measure), whose supports shrink to $\{(0,e)\}$ as $\delta \to 0$, and put
\[
    \psi_\delta(u,g) \; := \; \int \psi(u - a, hg) \, \chi_\delta(a,h) \, da \, dh   \qquad (u \in U_0).
\]
For small $\delta$ this is defined on $U_0 \times K^{\mathbb C}$, and $u - a \in U_1$ whenever $u \in U_0$ and $(a,h)$ lies in the support of $\chi_\delta$.
Since $(u,g) \mapsto (u-a, hg)$ is biholomorphic, $\psi_\delta$ is an average of plurisubharmonic functions and so is plurisubharmonic;
it is $\mathcal C^\infty$, as one sees by moving the translations onto $\chi_\delta$; and it is right $K$-invariant, since left and right
translations commute.
As $\psi$ is continuous, $\psi_\delta \to \psi$ uniformly on compact subsets of $U_0 \times K^{\mathbb C}$.

\medskip
Fix $\varepsilon > 0$ and put $\psi_{\delta,\varepsilon} := \psi_\delta + \varepsilon\rho$.
This is a $\mathcal C^\infty$, right $K$-invariant plurisubharmonic function on the tube $U_0 \times K^{\mathbb C}$, it is strictly
plurisubharmonic on each fiber, and it is an exhaustion on each fiber, because $\psi_\delta$ is bounded from below there.
By \cite[Th\'eor\`eme 3]{Loe85} the function $m_{\delta,\varepsilon}(u) := \inf_g \psi_{\delta,\varepsilon}(u,g)$ is plurisubharmonic on $U_0$.

\medskip
We claim that $m_{\delta,\varepsilon} \to m_\varepsilon(u) := \inf_g (\psi + \varepsilon\rho)(u,g)$ uniformly on $U_0$ as $\delta \to 0$.   
For $u \in \overline{U_0}$ let $g_u$ be a point at which $\psi(u,\cdot)$ attains its minimum.   
The points $(u,g_u)$ lie in the compact set $\{(u,g)\in\overline{U_0}\times K^{\mathbb C} : \psi(u,g)\le\max_{\overline{U_0}} m\}$, 
and so $m_\varepsilon(u) \le \psi(u,g_u) + \varepsilon\rho(g_u) \le c$ on $U_0$ for a suitable constant $c$.  
The set $A := \{(u,g) \in \overline{U_1}\times K^{\mathbb C} : \psi(u,g) \le c+1\}$ is compact.   
Let $B$ be a compact neighbourhood of $(0,e)$ containing the supports of all $\chi_\delta$ under consideration, and put 
$A_B := \{(u,g) \in \overline{U_0}\times K^{\mathbb C} : (u-a,hg) \in A \mbox{ for some } (a,h)\in B\}$, which is compact and contains 
$A \cap (\overline{U_0} \times K^{\mathbb C})$.     
If $u \in U_0$ and $(u,g) \notin A_B$, then $\psi(u-a,hg) > c+1$ for all $(a,h)$ in the support of $\chi_\delta$, 
so that $\psi_{\delta,\varepsilon}(u,g) \ge \psi_\delta(u,g) \ge c+1$, and also $(\psi+\varepsilon\rho)(u,g) > c+1$.    
On the other hand, if $\delta$ is so small that $|\psi_\delta - \psi| < 1$ on $A_B$, then $\psi_{\delta,\varepsilon}(u,g_u) < c+1$.   
Hence both infima are taken over $A_B$, so that $|m_{\delta,\varepsilon} - m_\varepsilon| \le \sup_{A_B} |\psi_\delta - \psi|$ on $U_0$, 
and the claim follows.    
A uniform limit of plurisubharmonic functions is plurisubharmonic, and so $m_\varepsilon$ is plurisubharmonic on $U_0$.

\medskip
Finally, $m_\varepsilon$ decreases as $\varepsilon \downarrow 0$, since $\rho \ge 0$, and $m \le m_\varepsilon \le \psi(u,g_u) + \varepsilon\rho(g_u)$,
where $g_u$ is a point at which $\psi(u,\cdot)$ attains its minimum.
Thus $m_\varepsilon \downarrow m$, and $m$, being the limit of a decreasing sequence of plurisubharmonic functions, is plurisubharmonic.
\end{proof}

\begin{thm}\label{ssmain}
Let $G$ be a connected complex Lie group, $\Gamma\subset G$ discrete with $G/\Gamma$ pseudoconvex, and 
$\mathcal A\lhd G$ a connected complex normal subgroup whose orbits in $G/\Gamma$ are relatively compact and positive dimensional 
(for instance the Hirschowitz annihilator of $G/\Gamma$ if it is normal and $G/\Gamma$ is not Stein). 
Put $L:={\rm cl}_G(\mathcal A\Gamma)$. 
Suppose the connected component of the identity $S:=L^0$ of $L$ is semisimple and put $I:=S^{\mathbb C}\Gamma$.
Then:
\begin{enumerate}[\rm(a)]
\item $I$ is a closed complex subgroup of $G$, $I/\Gamma$ is connected, $\dim I > 0$, and ${\mathcal O}(I/\Gamma)=\mathbb C$.
\item If $S$ has no compact simple factors, every continuous plurisubharmonic exhaustion $\varphi$ of $G/\Gamma$ is constant 
on the fibers of $G/\Gamma\to G/I$.
It therefore descends to a continuous plurisubharmonic exhaustion of $G/I$.   
In other words, $I$ satisfies (ii)$^*$ (Definition \ref{iistar}).
\item In general, $I$ satisfies (ii)$^*$ (Definition \ref{iistar}); 
here (MP) is used for $K:=S_c$, the product of the compact simple factors of $S$.
\end{enumerate}
\end{thm}

\begin{proof}
Let $\varphi$ be a continuous plurisubharmonic exhaustion of $G/\Gamma$ and $\tilde\varphi(g):=\varphi(g\Gamma)$.
Since $\mathcal A\lhd G$, the $\mathcal A$-orbit through $g\Gamma$ is $g\mathcal A\Gamma/\Gamma$.  
It is relatively compact, so, by Liouville's theorem along complex one-parameter subgroups of $\mathcal A$, $\varphi$ is constant on it, and by continuity
\begin{equation}\label{Linv0}
                           \varphi(g\ell\Gamma)  =  \varphi(g\Gamma)\qquad (g\in G,\ \ell\in L).
\end{equation}

\medskip 
(a) Since $\mathcal A\subset S$ and $S\Gamma$ is open, hence closed, in $L$, one has $L=S\Gamma$. 
Lemma \ref{sshull}, applied with $M=G$ and $P=L$, shows that $I=S^{\mathbb C}L=S^{\mathbb C}\Gamma$ 
is closed with $I^0=S^{\mathbb C}$.
Now let $f\in{\mathcal O}(I/\Gamma)$ and let $\tilde f$ be its lift to $I$.
Since the $\mathcal A$-orbits in the closed subset $I/\Gamma$ of $G/\Gamma$ are relatively compact, $f$ is bounded and hence constant 
on them (Liouville's theorem along complex one-parameter subgroups of $\mathcal A$), i.e., $\tilde f$ is right invariant under 
$\mathcal A\Gamma$, and by continuity under $L\supset S$.
But a holomorphic function annihilated by the left invariant vector fields of $\mathfrak s$ is also annihilated by those of $i\mathfrak s$, 
and so $\tilde f$ is right invariant under $S^{\mathbb C}\Gamma=I$, i.e., $f$ is constant.

\medskip
(b) By (\ref{Linv0}) the function $\tilde\varphi$ is right $S$-invariant, hence right $S^{\mathbb C}$-invariant by Proposition \ref{inv}, 
and, being also right $\Gamma$-invariant, it is right $I$-invariant.
Thus $\varphi$ descends to a continuous plurisubharmonic function $\psi$ on $G/I$, given locally by $\tilde\varphi\circ\sigma$ 
for local holomorphic sections $\sigma$ of $G\to G/I$, and $\{\psi\le c\}=\pi(\{\varphi\le c\})$ is compact, where $\pi:G/\Gamma\to G/I$.

\medskip
(c) Put $T:=S^{\mathbb C}=I^0$, a connected complex semisimple Lie group (see the proof of Proposition \ref{inv}), $K:=S_c$ and $T_{nc}:=S_{nc}^{\mathbb C}$, 
and let $p:T\to T/T_{nc}$ be as in Lemma \ref{proper}. 
By (\ref{Linv0}) and Proposition \ref{inv}, $\tilde\varphi$ is right invariant under $S$, $T_{nc}$ and $\Gamma$. 
Let $\sigma:U\to G$ be a holomorphic section of $G\to G/I$ over a coordinate ball $U\subset G/I$, and put 
\[
    \psi(u,k):=\tilde\varphi(\sigma(u)k)\qquad (u\in U,\ k\in K^{\mathbb C}).
\]
This is a continuous plurisubharmonic function on $U\times K^{\mathbb C}$, right invariant under $K\subset S$. 
Since $I=K^{\mathbb C}T_{nc}\Gamma$ (Lemma \ref{proper}(1)), $\inf_{k}\psi(u,k)$ is the infimum of $\varphi$ over the fiber $\sigma(u)I/\Gamma$, 
i.e., it equals $(m_I\varphi)(\sigma(u)I)$. 

\medskip
We next show that $\psi$ is an exhaustion relative to $U$ in the sense of (MP).   
Let $C\subset U$ be compact, $c\in\mathbb R$, and let $(u_n,k_n)\in C\times K^{\mathbb C}$ with $\psi(u_n,k_n)\le c$. 
Then $k_n\Gamma$ lies in the compact set $\sigma(C)^{-1}\cdot\{\varphi\le c\}$ and in the closed set $I/\Gamma$. 
As $T$ acts transitively on $I/\Gamma=T\Gamma/\Gamma$, there is a compact set $C'\subset T$ with $k_n\in C'\Gamma$, 
say $k_n=c_n\gamma_n$ with $c_n\in C'$ and $\gamma_n=c_n^{-1}k_n\in\Gamma\cap T$. 
Since $\Gamma\subset L$ normalizes $L^0=S$, Lemma \ref{proper}(2) gives $p(\gamma_n)\in p(K)$, so $p(k_n)$ lies in the compact set $p(C')p(K)$. 
By Lemma \ref{proper}(1) the sequence $(k_n)$ is then relatively compact in $K^{\mathbb C}$, which shows that the sublevel sets of $\psi$ over $C$ are compact.  

\medskip
By (MP), applied with $K^{\mathbb C}$ (which is closed by Lemma \ref{ssclosed}) and its compact real form $K$, 
the function $m_I\varphi$ is plurisubharmonic, and by Lemma \ref{fiberops} it is a continuous exhaustion of $G/I$. 
The same argument applies to a continuous plurisubharmonic function on $D\times G/\Gamma$ that is an exhaustion relative to $D$, 
with $U$ replaced by $D\times U$ (Remark \ref{starholds}), so $I$ satisfies (ii)$^*$.
\end{proof}

\section{Proof of the main result}   

\subsection{Formulation of the strategy of the proof}  
\begin{remark} \label{connected}  
There is a technical point that can arise in our construction,      
in that an intermediary fibration whose fiber is not connected might be involved.      
This is handled by a type of {\em Stein factorization for homogeneous fibrations}.      
Suppose $G$ is a connected Lie group that contains a closed subgroup $I$ containing a closed subgroup $H$.  
Let $\widetilde{I}$ be those connected components of $I$ that meet $H$.  
Then $\widetilde{I}$ is a closed subgroup of $G$ containing $H$ and  the fibration $G/H \to G/\widetilde{I}$ has connected fiber $\widetilde{I}/H$.    
\end{remark}

\noindent 
{\sc Strategy of the Proof:}

\medskip 
Assume $G/H$ is a pseudoconvex homogeneous manifold that is not Stein.  
In order to prove the theorem 
we construct a closed complex subgroup $I$ of $G$ containing $H$ with $\dim I > \dim H$ and $I/H$ connected,   
possibly with the aid of Remark \ref{connected}, such that    
\begin{enumerate}  
\item [] (i)  ${\mathcal O}(I/H) = \mathbb C$  and    
\item [] (ii) every continuous plurisubharmonic  
exhaustion function on $G/H$ induces a continuous plurisubharmonic exhaustion function on $G/I$.   
\end{enumerate}  
Then if $G/H \to G/J$ is the holomorphic reduction, $I$ is a subgroup of $J$ because of (i), $G/I$ is pseudoconvex because of (ii) 
and is either Stein (then $I=J$ and we are done) or not Stein and one applies the construction until 
one does reach the holomorphic reduction; see also the proof of Theorem \ref{mnthm} at the end of this section.  
In fact the subgroups $I$ will satisfy a stronger property (ii)$^*$, formulated below. 
By Lemma \ref{fiberops} it implies (ii), and unlike (ii) it passes from a closed subgroup $Q$ to $G$ (Lemma \ref{transfer}).

\begin{remark} 
The following complex homogeneous manifolds satisfy (ii) 
whenever they occur as the fiber of a homogeneous fibration $G/H \to G/I$.  
For (2) and (3) the basic tool is Kiselman's Minimum Principle \cite[Theorem 2.2]{Kis78}.  
\begin{enumerate}  
\item  compact complex homogeneous manifolds,       
\item  Cousin group towers, see \cite[Lemma 6.1 (1)]{GMO13},   
\item  the $\mathbb C^*$ power towers $R_L^{\mathbb C}\Gamma/\Gamma$ that occur in Step 1 of the proof of Proposition \ref{keyprop}, 
provided the exhaustion function is invariant under the totally real group $R_L$.        
\end{enumerate} 
\end{remark}

Several times in the proof we carry out the construction
not on $G/H$ itself but on a closed orbit $Q/H$ of a closed complex subgroup $Q$ of $G$;
this happens for $Q = G_1, G_2, G_3$ and for $Q = N$ below.
The construction then gives a subgroup $I$ with $H \subset I \subset Q$ that has
property (ii) for continuous plurisubharmonic exhaustions of $Q/H$.
What we need, however, is that every continuous plurisubharmonic exhaustion of $G/H$
induces one on $G/I$.
A function that is plurisubharmonic in each of two sets of variables separately need not be plurisubharmonic, 
so this does not follow by restricting to the fibers of $G/H \to G/Q$, 
and we have to keep track of how the induced function is constructed.
For this we use the following notation.

\medskip
For closed complex subgroups $H \subset I$ of $G$, a domain $D \subset \mathbb C^p$ (where $p \ge 0$,
and $D$ is a point if $p = 0$) and a function $f : D \times G/H \to \mathbb R$ we put
\[
     (m_I f)(u, gI) \; := \; \inf_{i \in I} f(u, giH) ,
\]
i.e., $m_I f$ is the infimum of $f(u,\cdot)$ over the fibers of $G/H \to G/I$.
Clearly $m_{I_1} \circ m_I = m_{I_1}$ whenever $H \subset I \subset I_1$.

\begin{definition} \label{iistar}
Let $H \subset I$ be closed complex subgroups of $G$, and let $D \subset \mathbb C^p$, $p \ge 0$, be a domain.
A function $\Phi$ on $D\times G/H$ is called an \emph{exhaustion relative to $D$} if for every compact set $C\subset D$
and every $c\in\mathbb R$ the set $\{(u,x)\in C\times G/H : \Phi(u,x)\le c\}$ is compact; for $p=0$ this is just an exhaustion.
We say that $I$ satisfies \textbf{(ii)$^*$ for $G/H$} if for every such $D$ and every continuous plurisubharmonic function $\Phi$
on $D \times G/H$ that is an exhaustion relative to $D$, the function $m_I\Phi$ is plurisubharmonic on $D \times G/I$.
\end{definition}

If $H \subset I \subset I_1$, where $I$ satisfies (ii)$^*$ for $G/H$ and $I_1$ satisfies (ii)$^*$ for $G/I$,
then $I_1$ satisfies (ii)$^*$ for $G/H$. 
This is because $m_{I_1}\Phi = m_{I_1}(m_I\Phi)$, and $m_I\Phi$ 
is again a continuous exhaustion relative to $D$ by the following lemma.

\begin{lemma} \label{fiberops}
Let $H \subset I$ be closed complex subgroups of $G$ and $f$ a continuous exhaustion function on $G/H$.
Then $m_I f$ is a continuous exhaustion function on $G/I$.
In particular, if $I$ satisfies (ii)$^*$ for $G/H$, then every continuous plurisubharmonic exhaustion $\varphi$ of $G/H$
induces the continuous plurisubharmonic exhaustion $m_I\varphi$ of $G/I$, i.e., (ii)$^*$ implies (ii). 
More generally, if $\Phi$ is a continuous exhaustion of $D \times G/H$ relative to $D$, then so is $m_I\Phi$ on $D\times G/I$.
\end{lemma}

\begin{proof}
Let $\pi : G/H \to G/I$ be the fibration.
Every fiber of $\pi$ is closed and $f$ is an exhaustion, so the infimum of $f$ over a fiber is attained
and one has $\{ m_I f \le c \} = \pi(\{ f \le c\})$, which is compact for every $c$.
In particular $m_I f$ is lower semicontinuous and an exhaustion.
Over a small open set $U \subset G/I$ with a holomorphic section $\sigma : U \to G$ one has
$(m_I f)(\sigma(v)I) = \inf_{i\in I} f(\sigma(v) i H)$, an infimum of continuous functions of $v$, and so $m_I f$
is also upper semicontinuous. 
The same argument, applied to ${\rm id}_C\times\pi : C\times G/H\to C\times G/I$ for a compact set $C\subset D$, 
gives the last statement:
the infimum of $\Phi(u,\cdot)$ over each fiber of $\pi$ is attained,
$\{(u,y)\in C\times G/I : (m_I\Phi)(u,y)\le c\}$ is the image of the compact set $\{(u,x)\in C\times G/H : \Phi(u,x)\le c\}$,
and over $D\times U$ one has $(m_I\Phi)(u,\sigma(v)I)=\inf_{i\in I}\Phi(u,\sigma(v)iH)$, 
an infimum of continuous functions of $(u,v)$.
\end{proof}

\begin{lemma}[Transfer Lemma] \label{transfer}
Let $Q$ be a closed complex subgroup of $G$ and $H \subset I \subset Q$ closed complex subgroups.
If $I$ satisfies (ii)$^*$ for $Q/H$, then $I$ satisfies (ii)$^*$ for $G/H$.
\end{lemma}

\begin{proof}
Let $D \subset \mathbb C^p$ be a domain and $\Phi$ a continuous plurisubharmonic function on $D \times G/H$ 
that is an exhaustion relative to $D$. 
Since plurisubharmonicity is a local property, it suffices to show that $m_I\Phi$ is plurisubharmonic on
$D \times W$ for $W$ running through an open covering of $G/I$.

\medskip
Fix $g_0 \in G$, let $V \subset \mathbb C^q$ be a ball and $s : V \to G$ a local holomorphic section of $G \to G/Q$,
i.e., $v \mapsto s(v)Q$ is a biholomorphic map of $V$ onto a neighbourhood of $g_0Q$ in $G/Q$.
For $K = H$ and $K = I$ the map
\[
    \alpha_K : V \times Q/K \; \longrightarrow \; G/K, \qquad (v, nK) \; \mapsto \; s(v)nK
\]
is injective (if $s(v_1)n_1K = s(v_2)n_2K$, then $s(v_1)Q = s(v_2)Q$, so $v_1 = v_2$ and $n_1K = n_2K$) and holomorphic,
and its image is the full preimage $W_K$ of $s(V)Q/Q$ in $G/K$; as source and target have the same dimension,
$\alpha_K$ is a biholomorphic map onto the open set $W_K$.
Moreover, $\alpha_H$ maps the fiber $\{v\} \times nI/H$ of $V \times Q/H \to V \times Q/I$ onto the fiber
$s(v)nI/H$ of $G/H \to G/I$.
Consequently,
\[
    (m_I\Phi)(u, \alpha_I(v,y)) \; = \; (m_I\Psi)(u,v,y) \qquad (u\in D, \ v \in V, \ y \in Q/I) ,
\]
where $\Psi(u,v,nH) := \Phi(u, s(v)nH)$ is a function on $(D\times V) \times Q/H$
and on the right the infimum is taken over the fibers of $Q/H \to Q/I$.

\medskip
Now $\Psi = \Phi \circ ({\rm id}_D \times \alpha_H)$ is plurisubharmonic.  
Moreover, $\Psi$ is a continuous exhaustion relative to
$D \times V$: if $C_1 \subset D$ and $C_2 \subset V$ are compact,
then $\{\Psi \le c\} \cap (C_1 \times C_2 \times Q/H)$ is the
image under $({\rm id}_D \times \alpha_H)^{-1}$ of the compact set
$\{\Phi \le c\} \cap (C_1 \times \pi_Q^{-1}(s(C_2)Q/Q))$, which lies in
the open set $D \times W_H$; here $\pi_Q : G/H \to G/Q$. 
Since $I$ satisfies (ii)$^*$ for $Q/H$, applied with the domain $D \times V \subset \mathbb C^{p+q}$, the function $m_I\Psi$ is
plurisubharmonic, and thus $m_I\Phi$ is plurisubharmonic on $D \times W_I$.
As $g_0$ was arbitrary, these sets cover $G/I$.
\end{proof}

The groups $G_1$, $G_2$, $G_3$ and $N$ that occur in the proof below need not be connected, 
and so we also need the following variant of Lemma \ref{transfer}, which allows us to work with their identity components.

\begin{lemma}[Component Lemma] \label{component}
Let $Q$ be a closed complex subgroup of $G$ containing $H$, put $H':=H\cap Q^0$, and let $I'$ be a closed complex subgroup of $Q^0$ 
containing $H'$ such that $hI'h^{-1}=I'$ for every $h\in H$. Put $I:=I'H$.
\begin{enumerate}[\rm(1)]
\item $I$ is a closed complex subgroup of $G$ with $H\subset I\subset Q$ and $I\cap Q^0=I'$, 
and $I'/H'\to I/H$, $i'H'\mapsto i'H$, is a biholomorphism. 
In particular, $\dim I-\dim H=\dim I'-\dim H'$, $I/H$ is connected if $I'/H'$ is, and ${\mathcal O}(I/H)={\mathcal O}(I'/H')$.
\item If $I'$ satisfies (ii)$^*$ for $Q^0/H'$, then $I$ satisfies (ii)$^*$ for $G/H$.
\end{enumerate}
\end{lemma}

\begin{proof}
(1) Since $H\subset Q$ normalizes $Q^0$ and $I'$, the set $I=I'H$ is a subgroup of $Q$, and $I\cap Q^0=I'(H\cap Q^0)=I'$. 
As $Q^0$ is open in $Q$, Lemma \ref{open}, applied with $M=Q$ and $U=Q^0$, shows that $I$ is closed; its identity component is $I'^0$, so it is complex. 
The map $i'H'\mapsto i'H$ is well defined, holomorphic and surjective, and it is injective because $i'^{-1}j'\in H\cap Q^0=H'$ whenever $i'H=j'H$. 
A bijective holomorphic map between complex manifolds is biholomorphic.

\medskip\noindent
(2) By Lemma \ref{transfer} it suffices to show that $I$ satisfies (ii)$^*$ for $Q/H$. 
Let $\Phi$ be a continuous plurisubharmonic function on $D\times Q/H$ that is an exhaustion relative to $D$. 
The connected components of $Q/H$ are the open and closed sets $qQ^0H/H$, $q\in Q$, and left translation by $q$ maps $Q^0H/H$ onto $qQ^0H/H$ 
and fibers of $Q/H\to Q/I$ onto fibers. 
So it suffices to consider $m_I\Phi$ on $D\times Q^0I/I$. 
Under the identifications $Q^0/H'\cong Q^0H/H$ and $Q^0/I'\cong Q^0I/I$, the fiber $q'I/H=q'I'H/H$, $q'\in Q^0$, of $Q/H\to Q/I$ 
corresponds to the fiber $q'I'/H'$ of $Q^0/H'\to Q^0/I'$. 
The restriction of $\Phi$ to the closed set $D\times Q^0H/H$ is a continuous plurisubharmonic function that is an exhaustion relative to $D$, 
and so $m_I\Phi=m_{I'}(\Phi|_{D\times Q^0/H'})$ is plurisubharmonic on $D\times Q^0/I'$ by assumption.
\end{proof}

\begin{remark} \label{starholds}
The subgroups $I$ constructed in the proof below satisfy (ii)$^*$, and not just (ii).
In each case the plurisubharmonicity of the induced function is proved in local trivializations $U \times I/H$
of $G/H \to G/I$ (or of intermediate fibrations), in which the coordinate $u$ of $U$ enters only as a parameter.
One of three arguments is used: 
the function is constant along the fibers (as in Theorem \ref{ssmain}(b)); 
Kiselman's minimum principle \cite[Theorem 2.2]{Kis78} applies, possibly in the form used in \cite[Lemma 6.1 (1)]{GMO13} for Cousin group towers; 
or (MP) applies. 
Each of these allows an arbitrary domain of parameters, and so $U$ may be replaced by $D \times U$.
The invariance properties needed for Kiselman's principle and for (MP),
namely (\ref{Linv0}), (\ref{Linv}) and Proposition \ref{inv}, are proved for an arbitrary continuous plurisubharmonic exhaustion
of the space in question. 
They therefore hold for each function $\Phi(u,\cdot)$ separately, 
and the sublevel sets of $\Phi$ over compact sets of parameters are compact, as (MP) requires.

\medskip
If a fiber is compact and connected, then by the maximum principle every plurisubharmonic function is constant on it, 
and so its infimum over the fiber is simply its value at any point of the fiber.

\medskip
For Cousin group towers we indicate why the proof of \cite[Lemma 6.1 (1)]{GMO13} gives (ii)$^*$. 
By the remark following Definition \ref{iistar}, and since each stage of a Cousin group tower is a homogeneous fibration, 
it suffices to treat a fibration $G/H\to G/I$ whose fiber $I/H$ is a Cousin group $C=\mathbb C^k/\Lambda$, $\Lambda\cong\mathbb Z^{k+l}$. 
Let $\Phi$ be a continuous plurisubharmonic function on $D\times G/H$ that is an exhaustion relative to $D$, and trivialize the bundle 
over a ball $B\subset G/I$ as $B\times C$. 
Then $\Phi$ lifts to a continuous plurisubharmonic function $\tilde\Phi(u,z,w)$ on $D\times B\times\mathbb C^k$ that is $\Lambda$-invariant in $w$. 
Let $\Lambda_{\mathbb R}\subset\mathbb C^k$ denote the real linear span of $\Lambda$, and put $V:=\Lambda_{\mathbb R}\cap i\Lambda_{\mathbb R}$, the largest complex subspace of $\Lambda_{\mathbb R}$. 
For fixed $(u,z)$ the function $w\mapsto\tilde\Phi(u,z,w)$ is bounded on every translate of $V$, since such a translate maps into a translate 
of the compact torus $\Lambda_{\mathbb R}/\Lambda$, and so it is invariant under $V$ by Liouville's theorem; 
as the image of $V$ is dense in $\Lambda_{\mathbb R}/\Lambda$, it is invariant under $\Lambda_{\mathbb R}$ by continuity. 
Thus $\tilde\Phi$ descends to $D\times B\times(\mathbb C^k/V)$ and is invariant under the real form $\Lambda_{\mathbb R}/V$ of $\mathbb C^k/V$, 
and Kiselman's minimum principle \cite[Theorem 2.2]{Kis78}, with $(u,z)$ as parameters, shows that 
$(u,z)\mapsto\inf_w\tilde\Phi(u,z,w)=(m_I\Phi)(u,z)$ is plurisubharmonic. 
Continuity and the exhaustion property relative to $D$ follow from Lemma \ref{fiberops}. 
This is the proof of \cite[Lemma 6.1 (1)]{GMO13}, with $(u,z)$ in place of $z$.
\end{remark}

\subsection{The proof itself}    

The main step in the proof of Theorem \ref{mnthm} is the following. 

\begin{proposition} \label{keyprop}
Let $G$ be a connected complex Lie group and $H$ a closed complex subgroup of $G$ such that $G/H$ is pseudoconvex
and not Stein.
Then there is a closed complex subgroup $I$ of $G$ containing $H$ with the following properties.
\begin{enumerate}
\item[(a)] $\dim I > \dim H$ and $I/H$ is connected.
\item[(b)] ${\mathcal O}(I/H) = \mathbb C$.
\item[(c)] $I$ satisfies (ii)$^*$ for $G/H$.
\item[(d)] $I$ is canonically attached to the pair $(G,H)$, i.e., every holomorphic automorphism $\alpha$
           of the complex Lie group $G$ with $\alpha(H) = H$ satisfies $\alpha(I) = I$.
\end{enumerate}
In particular, by Lemma \ref{fiberops} every continuous plurisubharmonic exhaustion $\varphi$ of $G/H$ induces the
continuous plurisubharmonic exhaustion $m_I\varphi$ of $G/I$, and so $G/I$ is pseudoconvex.
\end{proposition}

Property (d) is what makes the induction work. 
In Step 0 below an analogue of (d) shows that the subgroup constructed in a quotient of $G_3^0$ is normalized by $H$. 
In Case (C7) the Proposition is applied to the connected group $N^0$, and the resulting subgroup has to be normalized by $\Gamma$.

\begin{proof}[Proof of Proposition \ref{keyprop}]
The proof is rather long, and so we begin by listing the cases in the order in which they are treated, 
together with the subgroup $I$ obtained in each of them.  
\begin{enumerate}
\item[\rm (C1)] The fiber of $G_1/H\to G_1/G_2$ is Stein, where $G_1:=N_G(\mathcal A)$ and $G_2:=N_{G_1}(H^0)$. 
      Then the $\mathcal A$-orbits are flag manifolds and $I:=\mathcal AH$.
\item[] Otherwise we pass to a connected group $\widehat G$ with a discrete subgroup $\Gamma$ and a normal subgroup $\widehat{\mathcal A}$ 
      (Step~0); the group $I$ is recovered from a subgroup $\widehat I$ of $\widehat G$ by Lemma \ref{component}.
\item[\rm (C2)] $\widehat G$ is solvable: $\widehat I$ is the holomorphic reduction group; its fiber is a Cousin group tower.
\item[\rm (C3)] $\widehat G$ is semisimple: $\widehat I$ is the Remmert reduction group; its fiber is compact.
\item[] Otherwise $\widehat G$ is a mixed group, which we may take to be simply connected; we drop the hats, 
      and put $L:={\rm cl}_G(\mathcal A\Gamma)$ and $\widetilde L:=M_*\Gamma$ (Lemma \ref{Ltilde}).
\item[\rm (C4)] The $\mathcal A$-orbits are compact: $I:=\mathcal A\Gamma$.
\item[\rm (C5)] The radical $R_L$ of $L^0$ is trivial: $I:=(L^0)^{\mathbb C}\Gamma$, by Theorem \ref{ssmain}.
\item[\rm (C6)] $\widetilde H\Gamma/\Gamma$ is not Stein, where $\widetilde H\subset\widetilde L$ is given by Proposition \ref{Karl}: 
      $I$ is the holomorphic reduction group of $\widetilde H\Gamma/\Gamma$; its fiber is a Cousin group tower.
\item[\rm (C7)] $N\ne\widetilde L$, where $B:=(L\cap\widetilde H)^0$ and $N:=N_{\widetilde L}(B^{\mathbb C})^0\Gamma$: 
      $I:=I_N\Gamma$, where $I_N$ is given by the induction hypothesis for $(N^0,\Gamma\cap N^0)$. 
      (If $N/\Gamma$ is Stein we are in case (C4).)
\item[] Otherwise $B=R_L$, the group $R_L^{\mathbb C}\Gamma$ is closed, $R_L^{\mathbb C}\Gamma/\Gamma$ is a $\mathbb C^*$ power tower, 
      plurisubharmonic exhaustions descend to $G/R_L^{\mathbb C}\Gamma$ (Step 1), and $Y:=\widetilde L/R_L^{\mathbb C}\Gamma$ has 
      ${\mathcal O}(Y)=\mathbb C$ (Step 2).
\item[\rm (C8)] $\widetilde L_1^0$ is semisimple, where $\widetilde L_1:=\widetilde L/R_L^{\mathbb C}$: $Y$ is compact and $I:=\widetilde L$.
\item[\rm (C9)] Otherwise: $I$ is the preimage of $S_1^{\mathbb C}\Gamma_1$, where $S_1$ is the image of $L^0$ in $\widetilde L_1$; 
      here (MP) is used.
\end{enumerate}
The conditions that separate the cases depend only on $(G,H)$, and in each case the subgroup $I$ turns out to be canonical. 

\medskip
We argue by induction on $\dim G$.
If $\dim G = 0$, then $G/H$ is a point and so is Stein, and there is nothing to prove.
Assume that the Proposition holds for all pairs $(G_*, H_*)$ with $G_*$ a connected complex Lie group
of dimension less than $\dim G$.  
The Hirschowitz annihilator $\mathcal A$ of $G/H$ satisfies $\dim \mathcal A > \dim H$, because $G/H$ is not Stein \cite[Lemma 3.3]{GMO13}.    
Put $G_{1} := N_{G}(\mathcal A)$ and $G_{2} := N_{G_{1}}(H^0)$. 
These are closed complex subgroups of $G$ containing $H$ (\S\ref{Hirschowitz}), but they need not be connected. 
Because the $\mathcal A$-orbits are positive dimensional and relatively compact, the fiber $G_1/H$ of $G/H\to G/G_1$ is not Stein. 
By \S\ref{titssec} the base $G_1/G_2$ of $G_1/H\to G_1/G_2$ is an ${\rm Ad}(G_1)$-orbit in a Grassmannian, 
and hence it is equivariantly realized as an orbit in a projective space. 

\medskip\noindent
\emph{Case} (C1): \emph{the fiber $G_2/H$ of $G_1/H\to G_1/G_2$ is Stein} (e.g. if $G_1/H\to G_1/G_2$ is a covering). 
Since $\mathcal A$ is normal in $G_1$, the $\mathcal A$-orbits in $G_1/H$ are flag manifolds 
by Corollary \ref{flag2} (recall that connectedness of the group is not needed there). 
In particular $\mathcal A\cdot x_0$ is compact, so $I := \mathcal A H$ is closed, and $I/H=\mathcal A/(\mathcal A\cap H)$ is compact, connected 
and positive dimensional. 
Hence ${\mathcal O}(I/H)=\mathbb C$, every plurisubharmonic function is constant on the fibers of $G/H\to G/I$, and $I$ satisfies (ii)$^*$ 
(Remark \ref{starholds}). 
Finally, $I$ is canonical because $\mathcal A$ is.

\medskip\noindent
\emph{Step 0: reduction to $(\widehat G,\Gamma,\widehat{\mathcal A})$.} 
We may now assume that $G_2/H$ is not Stein.    
The connected components of $G_2/H$ are permuted transitively by left translations by elements of $G_2$, 
so the component $X_2:=G_2^0\cdot x_0\cong G_2^0/(H\cap G_2^0)$ is not Stein. 
It is pseudoconvex, being a closed submanifold of $G/H$. 
Let $\mathcal A_2\subset G_2^0$ be the Hirschowitz annihilator of $X_2$. 
Its orbits in $X_2$ are relatively compact and positive dimensional (\S\ref{Hirschowitz} and \cite[Lemma 3.3]{GMO13}). 
For $h\in H$, conjugation by $h$ is a holomorphic automorphism of $G_2^0$ preserving $H\cap G_2^0$, so it preserves $\mathcal A_2$ 
(\S\ref{Hirschowitz}). 
Hence $H\subset G_3:=N_{G_2}(\mathcal A_2)$, and $H^0\lhd G_3$ because $G_3\subset G_2=N_{G_1}(H^0)$.
Put
\[
   H':=H\cap G_3^0,\qquad \widehat{G} := G_3^0/H^0,\qquad \Gamma := H'/H^0,\qquad \widehat{\mathcal A} := \mathcal A_2/H^0 .
\]
Then $\widehat G$ is a connected complex Lie group, $\Gamma$ is a discrete subgroup and $\widehat{\mathcal A}\lhd\widehat G$ is connected. 
The manifold $\widehat G/\Gamma=G_3^0/H'$ is the connected component through $x_0$ of the closed submanifold $G_3/H$ of $G/H$, so it is pseudoconvex. 
The $\widehat{\mathcal A}$-orbits in it are leaves of the Levi foliation of $X_2$ lying in a closed subset, so they are relatively compact and positive 
dimensional; in particular $\widehat G/\Gamma$ is not Stein, and every plurisubharmonic function on it is constant on these orbits. 
The group $\widehat{\mathcal A}$ need not be the Hirschowitz annihilator of $\widehat G/\Gamma$, but we shall only use the properties just listed. 

\medskip
We shall construct a closed complex subgroup $\widehat I$ of $\widehat G$ containing $\Gamma$ such that 
\begin{enumerate}
\item[$(\widehat{\rm a})$] $\dim\widehat I>0$ and $\widehat I/\Gamma$ is connected;
\item[$(\widehat{\rm b})$] ${\mathcal O}(\widehat I/\Gamma)=\mathbb C$;
\item[$(\widehat{\rm c})$] $\widehat I$ satisfies (ii)$^*$ for $\widehat G/\Gamma$;
\item[$(\widehat{\rm d})$] $\alpha(\widehat I)=\widehat I$ for every holomorphic automorphism $\alpha$ of $\widehat G$ with $\alpha(\Gamma)=\Gamma$ 
and $\alpha(\widehat{\mathcal A})=\widehat{\mathcal A}$.
\end{enumerate}
Given $\widehat I$, let $I'\subset G_3^0$ be its preimage, a closed complex subgroup containing $H'$. 
For $h\in H$, conjugation by $h$ preserves $G_3^0$, $H^0$, $H'$ and $\mathcal A_2$, so it induces an automorphism of $\widehat G$ preserving $\Gamma$ and 
$\widehat{\mathcal A}$, and by $(\widehat{\rm d})$ it preserves $I'$. 
By the Component Lemma \ref{component}, applied with $Q=G_3$, the group $I:=I'H$ satisfies (a)--(c) for $G/H$. 
Finally, a holomorphic automorphism of $G$ stabilizing $H$ stabilizes $\mathcal A$, $G_1$, $G_2$, $X_2$, $\mathcal A_2$, $G_3$, $G_3^0$, $H'$ and $H^0$, 
and so it induces an automorphism of $\widehat G$ preserving $\Gamma$ and $\widehat{\mathcal A}$; by $(\widehat{\rm d})$ it preserves $I'$ and hence $I$, 
which proves (d). 

\medskip\noindent
\emph{Cases} (C2) \emph{and} (C3). 
If $\widehat G$ is solvable, let $\widehat I$ be the holomorphic reduction group of $\widehat G/\Gamma$, which is canonical. 
By \cite[Main Theorem and Theorem 8.5]{GMO13}, ${\mathcal O}(\widehat I/\Gamma)=\mathbb C$, so $\widehat I/\Gamma$ is connected. 
Since $\widehat G/\Gamma$ is not Stein and holomorphically separable solvmanifolds are Stein \cite{HO86}, $\widehat I\neq\Gamma$, and so $\dim\widehat I>0$. 
Moreover $\widehat I/\Gamma$ is a Cousin group tower \cite[Theorem 8.5]{GMO13}, and (ii)$^*$ holds by \cite[Lemma 6.1 (1)]{GMO13}, see Remark \ref{starholds}. 
If $\widehat G$ is semisimple, then $\widehat G/\Gamma$ is holomorphically convex \cite[Theorem 7.1]{GMO13}; 
let $\widehat I$ be its Remmert reduction group, which is canonical. 
The fiber $\widehat I/\Gamma$ is compact and connected, and positive dimensional because $\widehat G/\Gamma$ is not Stein, 
and (ii)$^*$ holds by Remark \ref{starholds}.

\medskip
It remains to treat the case where $\widehat{G}$ is a mixed group.   
Let $\pi:\widetilde G\to\widehat G$ be the universal covering. 
Replacing $\widehat G$, $\Gamma$ and $\widehat{\mathcal A}$ by $\widetilde G$, $\pi^{-1}(\Gamma)$ and $\pi^{-1}(\widehat{\mathcal A})^0$ does not change 
$\widehat G/\Gamma$ or the orbits of $\widehat{\mathcal A}$; 
automorphisms of $\widehat G$ preserving $\Gamma$ and $\widehat{\mathcal A}$ lift to automorphisms of $\widetilde G$ preserving the new groups, 
and a subgroup of $\widetilde G$ containing $\pi^{-1}(\Gamma)\supset\ker\pi$ is the preimage of its image. 
Thus we may also assume that $\widehat G$ is simply connected; 
the same device is used below whenever Proposition \ref{Karl} is applied to a closed subgroup of $\widehat G$.   
To simplify the notation we now drop the hats and write $G$, $\mathcal A$ and $I$ for $\widehat{G}$, $\widehat{\mathcal A}$ and $\widehat I$. 
In the rest of the proof, then, $G$ is a simply connected mixed group, $\Gamma\subset G$ is discrete, $G/\Gamma$ is pseudoconvex, 
and $\mathcal A\lhd G$ is a connected complex subgroup with relatively compact positive dimensional orbits. 
A subgroup is called canonical if it is invariant under all holomorphic automorphisms of $G$ preserving $\Gamma$ and $\mathcal A$.
 
\medskip  
Set $L := {\rm cl}_{G}(\mathcal A\cdot \Gamma)$.  
Since $\mathcal A\subset L^0$ and $L^0\Gamma$ is open, hence closed, in $L$, one has $L=L^0\Gamma$. 

\medskip\noindent
\emph{Case} (C4): \emph{the $\mathcal A$-orbits are compact}. 
Then $\mathcal A\Gamma/\Gamma$ is compact, hence closed, so $I:=\mathcal A\Gamma=L$ is a closed complex subgroup, 
$I/\Gamma$ is compact, connected and positive dimensional, ${\mathcal O}(I/\Gamma)=\mathbb C$, and $I$ satisfies (ii)$^*$ by Remark \ref{starholds}. 
Clearly $I$ is canonical. 
In the remaining cases the $\mathcal A$-orbits are not compact. 

\medskip
Let $\widetilde L:=M_*\Gamma$ be the closed complex subgroup given by Lemma \ref{Ltilde}.
Since the $\mathcal A$-orbits are relatively compact, $L/\Gamma$ is compact, and it follows that ${\mathcal O}(\widetilde{L}/\Gamma) = \mathbb C$.
To see this, let $f$ be a holomorphic function on $\widetilde L/\Gamma$. 
It is bounded, hence constant, on the $\mathcal A$-orbits, and so its lift $\tilde f$ is right invariant under $L$.
The closed subgroup of right translations leaving $\tilde f$ invariant thus contains the connected group $L^0$, and also $(L^0)^{\mathbb C}$.  
Following the construction in Lemma \ref{Ltilde}, we see that it contains 
$M_0$, ${\rm cl}_G(M_0\Gamma)$, $M_1$, \dots, and finally $\widetilde L$.

\medskip\noindent
\emph{Case} (C5): \emph{the radical $R_L$ of $L^0$ is trivial}. 
Since $\mathfrak{l}/\mathfrak{r}_L$ is semisimple and $L^0 \supset \mathcal{A}$ is not trivial, 
$L^0$ is then a non-trivial connected semisimple group. 
Theorem \ref{ssmain} gives the closed complex subgroup 
$I := (L^0)^{\mathbb C}\Gamma$ with $\mathcal{O}(I/\Gamma) = \mathbb C$ and $\dim I > 0$, 
and parts (b) and (c) of that theorem, applied with our normal subgroup $\mathcal A$, show that $I$ satisfies (ii)$^*$. 
As $L$ is canonical, so is $I$. 
In the rest of the proof we assume that $\dim R_L > 0$, a fact that will be used later.

\medskip   
We now turn to the remaining cases (C6)--(C9).    
Recall that $\widetilde L=M_*\Gamma$, where $M_*=\widetilde L^0$ is open in $\widetilde L$, and 
put $\Gamma_0:=\Gamma\cap M_*$, so that $\widetilde L/\Gamma=M_*/\Gamma_0$. 
By Lemma \ref{Ltilde} we have $R_L^{\mathbb C} \subset R_{\widetilde L}$, and so $\dim R_{M_*}>0$. 
We may therefore apply Proposition \ref{Karl} to the universal covering of $M_*$ and the preimage of $\Gamma_0$. 
Let $\widetilde H$ be the image in $M_*$ of the canonical group obtained there. 
Conjugation by $\gamma\in\Gamma$ is a holomorphic automorphism of $M_*$ preserving $\Gamma_0$; 
it lifts to the universal covering preserving the preimage of $\Gamma_0$, so it preserves $\widetilde H$. 
Thus all of $\Gamma$ normalizes $\widetilde H$, and $\widetilde H\Gamma$ is a group 
with $\widetilde H\Gamma\cap M_*=\widetilde H\Gamma_0$. 
The latter is closed in $M_*$: it is the image of the product of the canonical group with the preimage of $\Gamma_0$, 
which is closed by Proposition \ref{Karl} and contains the kernel of the covering map, i.e., is saturated. 
Since the covering map is open, $\widetilde H$ is open in $\widetilde H\Gamma_0$, hence it is the identity component of this closed group, 
and in particular $\widetilde H$ is closed. 
By Lemma \ref{open} (with $U=M_*$), $\widetilde H\Gamma$ is closed in $\widetilde L$, and there is a fibration   
\[   
       \widetilde{L}/\Gamma \; \longrightarrow \;  \widetilde{L}/ \widetilde{H}\cdot\Gamma     
\]    
with $ \widetilde{H}\cdot\Gamma$ closed in $ \widetilde{L}$, where $ \widetilde{H}$ is a connected, solvable, complex Lie group   
that contains the radical $R_{ \widetilde{L}}$ of $\widetilde{L}^0=M_*$ and is normalized by $\Gamma$.    
\medskip\noindent
\emph{Case} (C6): \emph{the fiber $\widetilde{H}\Gamma/\Gamma$ of this fibration is not Stein}. 
Then $\widetilde{H}\Gamma/\Gamma$ is a connected pseudoconvex solvmanifold that is not Stein.  
Let $I$ be its holomorphic reduction group, a closed complex subgroup of $\widetilde H\Gamma$ containing $\Gamma$; 
it is canonical because $\widetilde H$ is. 
By \cite[Theorem 8.5 and Definition 8.2]{GMO13} the fiber $I/\Gamma$ of the holomorphic reduction is a Cousin group tower, 
and it is positive dimensional because a holomorphically separable solvmanifold is Stein \cite{HO86}. 
Put $\Delta:=\Gamma\cap\widetilde H$ and let $J_0\subset\widetilde H$ be the holomorphic reduction group of $\widetilde H/\Delta\cong\widetilde H\Gamma/\Gamma$. 
Then $I\cap\widetilde H=J_0$ and $\Gamma\subset I$, so $I=J_0\Gamma$ and, as in Lemma \ref{component}(1), $I/\Gamma\cong J_0/\Delta$; 
thus the results of \cite{GMO13} for the solvmanifold $\widetilde H/\Delta$ apply to $I/\Gamma$. 
Consequently ${\mathcal O}(I/\Gamma) = \mathbb C$, and $I$ satisfies (ii)$^*$ for $\widetilde H\Gamma/\Gamma$ by \cite[Lemma 6.1 (1)]{GMO13} along the tower 
(Remark \ref{starholds}), and hence for $G/\Gamma$ by Lemma \ref{transfer}, applied with $Q=\widetilde{H}\Gamma$. 
In the remaining cases the fiber $\widetilde{H}\Gamma/\Gamma$ is Stein. 

\medskip  
Now in order to finish the proof that $\widetilde{L}/\Gamma$ satisfies (ii) we have to analyze the structure of (at least part of)    
the intersection of the compact orbit $L/\Gamma$ with the Stein orbit $\widetilde{H}/\widetilde{H}\cap\Gamma$.      
Consider the connected, real, solvable group $B := (L \cap \widetilde{H})^0$.  
Since $R_L \subset R_{\widetilde{L}}\subset\widetilde{H}$ by Lemma \ref{Ltilde} and $R_L \subset L$, it follows that 
$R_L \subset B$ and thus $\dim B >0$.   
Let $B^{\mathbb C}$ be its complexification.      
We have $\Gamma \subset N_{\widetilde{L}}(B^{\mathbb C})$,   
since $\Gamma \subset L$ and $\Gamma$ normalizes $\widetilde{H}$ by Proposition \ref{Karl},   
and thus we may consider the fibration $\widetilde{L}/\Gamma \to \widetilde{L}/N_{\widetilde{L}}(B^{\mathbb C})$.   

\medskip
Now put $N' := N_{\widetilde{L}}(B^{\mathbb C})$, a closed complex subgroup of $\widetilde L$ containing $\Gamma$, 
and let $N:=N'^0\Gamma$ be the union of the connected components of $N'$ that meet $\Gamma$ (Remark \ref{connected}).
The base $\widetilde L/N'$ of the fibration $\widetilde L/\Gamma \to \widetilde L/N'$ is the ${\rm Ad}(\widetilde L)$-orbit of $\mathfrak b^{\mathbb C}$ 
in a Grassmannian, and so is equivariantly realized as an orbit in a projective space.
The connected components of $N'/\Gamma$ are the sets $nN/\Gamma$, $n\in N'$, and left translation by $n$ maps $N/\Gamma$ biholomorphically onto $nN/\Gamma$. 
So if $N/\Gamma$ were Stein, $N'/\Gamma$ would be holomorphically separable, and the $\mathcal A$-orbits would be compact by Corollary \ref{flag2},
applied with $G$, $H$, $I$, $J$ replaced by $\widetilde L$, $\Gamma$, $N'$, $\mathcal A$ 
(recall that $\mathcal A \lhd \widetilde L$); this is Case (C4). 
Hence $N/\Gamma$ is not Stein.

\medskip
\noindent
\emph{Case} (C7): \emph{$N \neq \widetilde L$}.
Then $N^0 \neq \widetilde L^0 = M_*$, since otherwise $B^{\mathbb C}$ would be normalized by both $M_*$ and $\Gamma$,
and thus $\dim N^0 < \dim \widetilde L \le \dim G$.
Set $\Gamma_N := \Gamma \cap N^0$.
Since $N = N^0\Gamma$, the map $N^0/\Gamma_N \to N/\Gamma$ is an $N^0$-equivariant biholomorphism, and so
$N^0/\Gamma_N$ is pseudoconvex (being a closed submanifold of the pseudoconvex manifold $\widetilde L/\Gamma$) and not Stein.
By the induction hypothesis there is a closed complex subgroup $I_N$ of $N^0$ containing $\Gamma_N$
satisfying (a)--(d) for the pair $(N^0, \Gamma_N)$.
For $\gamma \in \Gamma$ conjugation by $\gamma$ is a holomorphic automorphism of $N^0$ that stabilizes $\Gamma_N$,
and so by (d) it stabilizes $I_N$.
Thus $I := I_N\Gamma$ is a subgroup of $N$, and $I \cap N^0 = I_N$, because $\Gamma \cap N^0 = \Gamma_N \subset I_N$.
By Lemma \ref{open}, applied with $U = N^0$, the group $I$ is closed.
As above, the maps $N^0/\Gamma_N \to N/\Gamma$ and $N^0/I_N \to N/I$ are $N^0$-equivariant
biholomorphisms that are compatible with the fibrations $N^0/\Gamma_N \to N^0/I_N$ and $N/\Gamma \to N/I$.
Hence $I/\Gamma \cong I_N/\Gamma_N$ is connected, positive dimensional and has only constant holomorphic functions,
and $I$ satisfies (ii)$^*$ for $N/\Gamma$.
By Lemma \ref{transfer}, applied with $Q = N$, it also satisfies (ii)$^*$ for $G/\Gamma$.
Finally, a holomorphic automorphism of $G$ stabilizing $\Gamma$ and $\mathcal A$ stabilizes $L$, $\widetilde L$, $\widetilde H$, $B$, $N'$ and $N$,
and so it restricts to a holomorphic automorphism of $N^0$ stabilizing $\Gamma_N$; by (d) for $(N^0,\Gamma_N)$
it stabilizes $I_N$, and so also $I$.
This settles Case (C7).

\medskip 
In the remaining cases $N = \widetilde L$, i.e., $B^{\mathbb C}$ is a connected, normal, solvable subgroup of ${\widetilde{L}}$.  
Then $B^{\mathbb C} \cap L$ is normal in $L$, which implies that $B \subset R_L$.  
But $R_L \subset B$ and so $B=R_L$.      
We claim that the orbit $R_L\Gamma/\Gamma$ of $R_L$ through the base point is compact.    
As noted above, $R_L \subset R_{\widetilde{L}}$ by Lemma \ref{Ltilde}.  
Since $L\cap\widetilde H\Gamma = (L\cap\widetilde H)\Gamma$ is closed, $L\cap\widetilde H$ is open in it, and the identity component 
of $L\cap\widetilde H$ is $B=R_L$, we have $\mathrm{cl}(R_L\Gamma)^0 \subset R_L$. 
So $R_L$ is open in $\mathrm{cl}(R_L\Gamma)$ and, as in the proof of Lemma \ref{Ltilde}, $R_L\Gamma$ is closed. 
It follows that the orbit $R_L\Gamma/\Gamma$ is closed, and since it lies in the compact set $L/\Gamma$, it is compact.   
If ${\mathfrak r}_L \cap i {\mathfrak r}_L \not= (0)$, then the connected complex Lie group corresponding to 
the complex ideal ${\mathfrak r}_L \cap i {\mathfrak r}_L $ has positive dimensional orbits in the compact orbit $R_L\Gamma/\Gamma$, 
which is contained in the Stein manifold $\widetilde H\Gamma/\Gamma$.  
By the maximum principle this is impossible, since the $\widetilde{H}$-orbits are now Stein.  
Therefore ${\mathfrak r}_L \cap i {\mathfrak r}_L = (0)$, and the orbit $R_L\Gamma/\Gamma$ is totally real.  
Let $R_L^{\mathbb C}$ be the complexification of $R_L$, with 
Lie algebra ${\mathfrak r}_L^{\mathbb C} :=  {\mathfrak r}_L \oplus i {\mathfrak r}_L$.   
Since $R_L\subset\widetilde H$ and $\widetilde H$ is complex, $R_L^{\mathbb C}\subset\widetilde H$. 
We now show that $R_L^{\mathbb C}\Gamma$ is closed in $\widetilde L$.
The group $\widetilde H$ is closed (Proposition \ref{Karl}) and normal in the closed group $\widetilde H\Gamma$, 
and $\widetilde H\Gamma/\widetilde H$ is countable, so $\widetilde H$ is open in $\widetilde H\Gamma$.
Since $\Gamma$ normalizes $R_L^{\mathbb C}$, the set $R_L^{\mathbb C}\Gamma$ is a subgroup of $\widetilde H\Gamma$, 
and $R_L^{\mathbb C}\Gamma\cap\widetilde H=R_L^{\mathbb C}\Delta$ with $\Delta:=\Gamma\cap\widetilde H$.
By Lemma \ref{open} it suffices to show that $R_L^{\mathbb C}\Delta$ is closed in $\widetilde H$.
Let $E\to\widetilde H$ be the universal covering, $\Delta_E$ the preimage of $\Delta$ 
and $W$ the identity component of the preimage of $R_L$.
Then $E$ is a simply connected solvable complex Lie group, $E/\Delta_E=\widetilde H/\Delta$ is Stein, $\mathfrak w\cap i\mathfrak w=0$, 
$W^{\mathbb C}\lhd E$ because $R_L^{\mathbb C}\lhd\widetilde L$ (Lemma \ref{Ltilde}), $\Delta_E$ normalizes $W$, 
and $W\cap\Delta_E$ is cocompact in $W$ because $W/(W\cap\Delta_E)\cong R_L/(R_L\cap\Gamma)\cong R_L\Gamma/\Gamma$ is compact.
By Theorem \ref{solv}, $W^{\mathbb C}\Delta_E$ is closed, and so is its image $R_L^{\mathbb C}\Delta$, the preimage being saturated.
Hence $R_L^{\mathbb C}\Gamma$ is closed, and the $R_L^{\mathbb C}$-orbits are closed in $\widetilde{L}/\Gamma$.
Moreover, Theorem \ref{solv} shows that these orbits are $\mathbb C^*$ power towers, each of whose stages is Stein.

\medskip
We now show that $R_L^{\mathbb C}\Gamma$ satisfies (ii)$^*$ for $G/\Gamma$ (Step 1), and then analyze the fiber of $G/R_L^{\mathbb C}\Gamma\to G/\widetilde L$ (Step 2).
Let $\varphi$ be a continuous plurisubharmonic exhaustion function on $G/\Gamma$.
It is constant on each $\mathcal A$-orbit, by Liouville's theorem applied along complex one-parameter subgroups of $\mathcal A$. 
Since $\mathcal A\lhd G$, the $\mathcal A$-orbit through $g\Gamma$ is $g\mathcal A\Gamma/\Gamma$, and its closure is $gL/\Gamma$.
By continuity we get
\begin{equation}\label{Linv}
                                \varphi(g\ell\Gamma)\;  =  \;\varphi(g\Gamma)\qquad\mbox{for all } g\in G,\ \ell\in L .
\end{equation}

\medskip
{\it Step 1: Dividing out by $R_L^{\mathbb C}$.} 
Put $\mathfrak r_0:=0$, and for $i\ge1$ let $\mathfrak r_i\subset\mathfrak r_L$ be the preimage of the center of the nilradical of $\mathfrak r_L/\mathfrak r_{i-1}$. 
Each $\mathfrak r_i$ is a characteristic ideal of $\mathfrak r_L$, and $\mathfrak r_i/\mathfrak r_{i-1}$ is abelian and non-zero as long as $\mathfrak r_{i-1}\neq\mathfrak r_L$, 
so that $\mathfrak r_m=\mathfrak r_L$ for some $m$. 
Let $R_i$ be the connected subgroup of $R_L$ with Lie algebra $\mathfrak r_i$ and put $Z_i:=R_i^{\mathbb C}$, whose Lie algebra is 
$\mathfrak z_i=\mathfrak r_i\oplus i\mathfrak r_i$, because $\mathfrak r_L\cap i\mathfrak r_L=0$. 
Thus $Z_0=\{e\}$, $Z_m=R_L^{\mathbb C}$, $Z_{i-1}\lhd Z_i$, and $Z_i/Z_{i-1}$ is abelian. 
Since $L$ normalizes $R_L$, it normalizes every $R_i$ and every $Z_i$; in particular so does $\Gamma$. 
Applied to $(E,\Delta_E,W)$, the inductive proof of Theorem \ref{solv} runs exactly through the Lie algebras $\mathfrak z_i$, 
and it shows that at each stage the corresponding connected subgroup of $E$ times $\Delta_E$ is closed. 
Passing to the images in $\widetilde H$ and using Lemma \ref{open} as above, we see that every $Z_i\Gamma$ is closed in $\widetilde L$, hence in $G$.

\medskip\noindent
\emph{Claim.} Let $1\le i\le m$, let $D\subset\mathbb C^p$ be a domain, and let $\Phi$ be a continuous plurisubharmonic function on 
$D\times G/Z_{i-1}\Gamma$ that is an exhaustion relative to $D$ and satisfies 
\begin{equation}\label{Linvi}
\Phi(u,g\ell Z_{i-1}\Gamma)=\Phi(u,gZ_{i-1}\Gamma)\qquad(u\in D,\ g\in G,\ \ell\in L).
\end{equation}
Then $m_{Z_i\Gamma}\Phi$ is a continuous plurisubharmonic function on $D\times G/Z_i\Gamma$ that is an exhaustion relative to $D$ 
and again satisfies (\ref{Linvi}), with $Z_i$ in place of $Z_{i-1}$.

\medskip\noindent
\emph{Proof of the Claim.} 
Choose a real subspace $\mathfrak v_{\mathbb R}\subset\mathfrak r_i$ complementary to $\mathfrak r_{i-1}$ and put 
$\mathfrak v:=\mathfrak v_{\mathbb R}\oplus i\mathfrak v_{\mathbb R}$, a complex subspace of $\mathfrak z_i$ complementary to $\mathfrak z_{i-1}$. 
Since $Z_i/Z_{i-1}$ is a connected abelian Lie group with Lie algebra $\mathfrak z_i/\mathfrak z_{i-1}\cong\mathfrak v$, the map 
$\mathfrak v\to Z_i/Z_{i-1}$, $w\mapsto\exp(w)Z_{i-1}$, is a surjective homomorphism; in particular 
$\exp(w+x)Z_{i-1}=\exp(w)\exp(x)Z_{i-1}$ for $w,x\in\mathfrak v$. 
Let $\sigma:U\to G$ be a holomorphic section of $G\to G/Z_i\Gamma$ over a coordinate ball $U$. 
Because $\Gamma$ normalizes $Z_{i-1}$, the fiber of $G/Z_{i-1}\Gamma\to G/Z_i\Gamma$ over $\sigma(v)Z_i\Gamma$ is 
$\{\sigma(v)\exp(w)Z_{i-1}\Gamma : w\in\mathfrak v\}$. 
Put 
\[
     F(u,v,w):=\Phi\bigl(u,\sigma(v)\exp(w)Z_{i-1}\Gamma\bigr)\qquad(u\in D,\ v\in U,\ w\in\mathfrak v).
\]
This is a continuous plurisubharmonic function. 
If $x\in\mathfrak v_{\mathbb R}$, then $\exp(x)\in R_i\subset L$, and (\ref{Linvi}) with $g=\sigma(v)\exp(w)$ and $\ell=\exp(x)$ gives $F(u,v,w+x)=F(u,v,w)$. 
So, writing $w=a+ib$ with $a,b\in\mathfrak v_{\mathbb R}$, the function $F$ does not depend on $a$, and by Kiselman's minimum principle \cite[Theorem 2.2]{Kis78} 
the function $(u,v)\mapsto\inf_{b\in\mathfrak v_{\mathbb R}}F(u,v,ib)$ is plurisubharmonic on $D\times U$. 
This infimum is $(m_{Z_i\Gamma}\Phi)(u,\sigma(v)Z_i\Gamma)$. 
By Lemma \ref{fiberops}, $m_{Z_i\Gamma}\Phi$ is continuous and an exhaustion relative to $D$. 
Finally, if $\ell\in L$ and $z\in Z_i$, then $\ell z=z'\ell$ with $z':=\ell z\ell^{-1}\in Z_i$, and (\ref{Linvi}) gives 
$\Phi(u,g\ell zZ_{i-1}\Gamma)=\Phi(u,gz'Z_{i-1}\Gamma)$; taking the infimum over $z\in Z_i$ proves the invariance of $m_{Z_i\Gamma}\Phi$. 
This proves the claim.

\medskip
Now let $\Phi$ be a continuous plurisubharmonic function on $D\times G/\Gamma$ that is an exhaustion relative to $D$. 
Each $\Phi(u,\cdot)$ is a continuous plurisubharmonic function on $G/\Gamma$, so it satisfies (\ref{Linv}), i.e., $\Phi$ satisfies (\ref{Linvi}) for $i=1$. 
Applying the claim for $i=1,\dots,m$ and using $m_{Z_m\Gamma}=m_{Z_m\Gamma}\circ m_{Z_{m-1}\Gamma}\circ\dots\circ m_{Z_1\Gamma}$, 
we see that $R_L^{\mathbb C}\Gamma=Z_m\Gamma$ satisfies (ii)$^*$ for $G/\Gamma$. 
In particular, for a continuous plurisubharmonic exhaustion $\varphi$ of $G/\Gamma$, the function $\varphi_R:=m_{R_L^{\mathbb C}\Gamma}\varphi$ 
is a continuous plurisubharmonic exhaustion of $G/R_L^{\mathbb C}\Gamma$ satisfying (\ref{Linv}).

\medskip 
{\it Step 2: The fiber $Y$.} 
By Lemma \ref{Ltilde}, $R_L^{\mathbb C}\lhd\widetilde L$.
Put $\widetilde L_1:=\widetilde L/R_L^{\mathbb C}$, and let $\Gamma_1$ and $L_1$ be the images of $\Gamma$ and $L$.
Since $R_L^{\mathbb C}\Gamma$ is closed, $\Gamma_1$ is a countable closed subgroup of $\widetilde L_1$, and hence it is discrete.
The fiber $Y:=\widetilde L/R_L^{\mathbb C}\Gamma\cong\widetilde L_1/\Gamma_1$ of $G/R_L^{\mathbb C}\Gamma\to G/\widetilde L$ is closed.
Hence $\varphi_R|_Y$ is a continuous plurisubharmonic exhaustion and $Y$ is pseudoconvex.
Also ${\mathcal O}(Y)=\mathbb C$, since ${\mathcal O}(\widetilde L/\Gamma)=\mathbb C$, 
and $L_1/\Gamma_1$ is compact, being the image of $L/\Gamma$.
The group $(L\cap R_L^{\mathbb C})^0$ is connected, solvable and normal in $L$, so it equals $R_L$.  
The group $L_1$ is closed in $\widetilde L_1$, since $L_1\supset\Gamma_1$ is the preimage of the compact set $L_1/\Gamma_1$.
Being the continuous image of a connected Lie group,
the image of $L^0$ is $\sigma$-compact and has countable index in $L_1$, because $L=L^0\Gamma$ and $\Gamma$ is countable. 
Since $L_1$ is locally compact, Baire's theorem shows 
that the image of $L^0$ has nonempty interior and is therefore open in $L_1$. 
Hence the identity component of $L_1$ is the image of $L^0$ and has Lie algebra 
$\mathfrak l/(\mathfrak l\cap\mathfrak r_L^{\mathbb C})=\mathfrak l/\mathfrak r_L$, which is semisimple.
Thus dividing once by $R_L^{\mathbb C}$ removes the radical of $L$, and the triple $(\widetilde L_1,\Gamma_1,L_1)$ 
is in the situation where $L_1^0$ is semisimple.

\medskip 
{\it Step 3: Conclusion.} 
There are two cases.

\medskip\noindent
\emph{Case} (C8): \emph{$\widetilde L_1^0$ is semisimple}. 
Then $Y$ is holomorphically convex by \cite[Theorem 7.1]{GMO13}, applied to the connected semisimple group $M_*/R_L^{\mathbb C}=\widetilde L_1^0$, 
which acts transitively on $Y$. Since ${\mathcal O}(Y)=\mathbb C$, it follows that $Y$ is compact.
The homogeneous fibration $G/R_L^{\mathbb C}\Gamma\to G/\widetilde L$ is a holomorphic fiber bundle with compact fiber $Y$.
Since $Y$ is compact and connected, $\varphi_R$ is constant on the fibers (Remark \ref{starholds}). 
So $\varphi_2:=m_{\widetilde L}\varphi_R=m_{\widetilde L}\varphi$ is given in a local trivialization $U\times Y$ by $u\mapsto\varphi_R(u,y)$ for any $y\in Y$; 
it is plurisubharmonic, and by Lemma \ref{fiberops} it is a continuous exhaustion of $G/\widetilde L$. 
In this case $I:=\widetilde L$ satisfies (a)--(c), and it is canonical.

\medskip
\noindent
\emph{Case} (C9): \emph{$\widetilde L_1^0$ is not semisimple}. 
Let $S_1:=L_1^0$, which is semisimple, and $I_1:=S_1^{\mathbb C}\Gamma_1$.
By Lemma \ref{sshull} with $M=\widetilde L_1$ and $P=L_1$, the group $I_1$ is 
a closed complex subgroup of $\widetilde L_1$ with $I_1^0=S_1^{\mathbb C}$.
Let $I\subset\widetilde L$ be its preimage, a closed complex subgroup containing $R_L^{\mathbb C}\Gamma$ with $I/\Gamma$ connected.

\medskip 
It remains to show that $I$ satisfies (b) and (c) of Proposition \ref{keyprop}; (a) is clear, since $I\supset R_L^{\mathbb C}$.  
For (b) let $f\in{\mathcal O}(I/\Gamma)$.
As in Theorem \ref{ssmain}(a), the function $f$ is bounded and hence constant on the relatively compact complex $\mathcal A$-orbits 
in the closed subset $I/\Gamma\subset G/\Gamma$. 
Its lift $\tilde f$ is therefore right invariant under $\mathcal A\Gamma$.   
Hence $\tilde f$ is right invariant under $L$, and, by holomorphy, under $R_L^{\mathbb C}$,  
so that $f$ descends to $I_1/\Gamma_1$. 
There the descended function is right invariant under $L_1\supset S_1$, hence under $S_1^{\mathbb C}\Gamma_1=I_1$, so $f$ is constant.

\medskip
For (c), observe first that the function $\tilde\varphi_R$ on $G$ is right invariant under $L\,R_L^{\mathbb C}$.
Let $S$ be a Levi factor of $L^0$, so that $S$ maps onto $S_1$ with discrete kernel.
By Proposition \ref{inv} with $M=G$, $\tilde\varphi_R$ is right invariant under $S_{nc}^{\mathbb C}$.  
Together with $R_L^{\mathbb C}$ and $\Gamma$ it is therefore right invariant under the preimage of $(S_1)_{nc}^{\mathbb C}\Gamma_1$.
If $S_1$ has no compact simple factors, $\tilde\varphi_R$ is right invariant under $I$, and as in Theorem \ref{ssmain}(b),  
it induces a continuous plurisubharmonic exhaustion of $G/I$.
In general, let $q:\widetilde L\to\widetilde L_1$ be the quotient map, $K:=S_c$, and $K_1:=q(K)=(S_1)_c$. 
Since $q(S^{\mathbb C})=S_1^{\mathbb C}$ and $\ker q=R_L^{\mathbb C}$, one has $I=S^{\mathbb C}R_L^{\mathbb C}\Gamma$, 
and by Lemma \ref{proper}(1), applied with $T=S^{\mathbb C}$, $S^{\mathbb C}=S_{nc}^{\mathbb C}K^{\mathbb C}$ with commuting factors. 
Let $\sigma:U\to G$ be a holomorphic section of $G\to G/I$ over a coordinate ball $U$ and put $\psi(u,k):=\tilde\varphi_R(\sigma(u)k)$ 
for $(u,k)\in U\times K^{\mathbb C}$. 
As $\tilde\varphi_R$ is right invariant under $L$, $R_L^{\mathbb C}$, $\Gamma$ and $S_{nc}^{\mathbb C}$, the function $\psi$ is continuous, plurisubharmonic 
and right $K$-invariant, and $\inf_k\psi(u,k)$ is the infimum of $\varphi_R$ over the fiber of $G/R_L^{\mathbb C}\Gamma\to G/I$ through $\sigma(u)$. 

\medskip
As in the proof of Theorem \ref{ssmain}(c), $\psi$ is an exhaustion relative to $U$. 
To see this, let $C\subset U$ be compact, $c\in\mathbb R$ and $(u_n,k_n)\in C\times K^{\mathbb C}$ with $\psi(u_n,k_n)\le c$. 
There is a compact set $C'\subset S^{\mathbb C}$ with $k_n\in C'R_L^{\mathbb C}\Gamma$, 
since $S^{\mathbb C}$ acts transitively on $I/R_L^{\mathbb C}\Gamma$. 
Then $q(k_n)=q(c_n)\gamma_n'$ with $c_n\in C'$ and $\gamma_n'\in\Gamma_1\cap S_1^{\mathbb C}$, which normalizes $S_1=L_1^0$ 
because $\Gamma_1\subset L_1$. 
With $p_1:S_1^{\mathbb C}\to S_1^{\mathbb C}/(S_1)_{nc}^{\mathbb C}$, Lemma \ref{proper}(2) gives $p_1(q(k_n))\in p_1(q(C'))p_1(K_1)$, a compact set. 
The homomorphism $q|_{K^{\mathbb C}}:K^{\mathbb C}\to K_1^{\mathbb C}$ is surjective with kernel $K^{\mathbb C}\cap R_L^{\mathbb C}$, which is discrete 
(a solvable ideal of $\mathfrak k^{\mathbb C}$ is trivial) and normal, hence central and finite; together with Lemma \ref{proper}(1) 
it follows that $p_1\circ q|_{K^{\mathbb C}}$ is proper. 
Consequently $(k_n)$ is relatively compact in $K^{\mathbb C}$.  

\medskip
By (MP) the fiberwise infimum is plurisubharmonic, and it is attained, since the fibers are closed and $\varphi_R$ is an exhaustion. 
By Lemma \ref{fiberops} it is a continuous plurisubharmonic exhaustion of $G/I$, and so $I$ satisfies (ii)$^*$. 
Finally, every group that enters the definition of $I$ (the subgroup $\mathcal A$,
normalizers, closures, complexifications, radicals and the centers of
their nilradicals, the group $\widetilde H$ of Proposition \ref{Karl}, the image
$L_1^0$ of $L^0$, and holomorphic reduction groups) is canonically attached to the
triple $(G,\Gamma,\mathcal A)$, and so $I$ is canonical in the sense fixed above.  A Levi factor $S$ of $L^0$,
which is not canonical, is used only in the proof of (ii), not in the definition of $I$. 
This completes the proof of the Proposition.
\end{proof}

\begin{proof}[Proof of Theorem \ref{mnthm}]
Let $G/H$ be pseudoconvex with holomorphic reduction $G/H \to G/J$.
Put $J_0 := H$.
If $G/J_k$ is pseudoconvex and not Stein, let $J_{k+1}$ be the group $I$ given by Proposition \ref{keyprop} for the pair $(G, J_k)$.
Then $G/J_{k+1}$ is again pseudoconvex and $\dim J_{k+1} > \dim J_k$, and so after finitely many steps one arrives at
a group $J_m$ for which $G/J_m$ is Stein.

\medskip
We claim that $J_m = J$.
First, ${\mathcal O}(J_k/H) = \mathbb C$ for every $k$.
This is clear for $k=0$, and if it holds for $k$, then a holomorphic function on $J_{k+1}/H$ is constant on the fibers
$J_k/H$ of $J_{k+1}/H \to J_{k+1}/J_k$, so it is the pull back of a holomorphic function on $J_{k+1}/J_k$, which is constant by (b).
As a consequence, every $f \in {\mathcal O}(G/H)$ is constant on the fibers $gJ_m/H$ of $G/H \to G/J_m$, i.e., $J_m \subset J$.
On the other hand, the holomorphic functions on the Stein manifold $G/J_m$ separate its points, and their pull backs to $G/H$
lie in ${\mathcal O}(G/H)$ and are constant on the fibers of $G/H \to G/J$, so that $J \subset J_m$.
Thus $J = J_m$, the base $G/J$ is Stein and ${\mathcal O}(J/H) = \mathbb C$.
(If ${\mathcal O}(G/H) = \mathbb C$, this says that $J = G$, which is of course also immediate.)
\end{proof}

\end{document}